\documentclass[preprint, a4paper, 11pt, oneside, onecolumn]{elsarticle}
\usepackage{amssymb, amsthm,amsmath}
\usepackage{graphicx, subcaption, placeins}
\usepackage{float}    
\usepackage{amsthm}
\usepackage{amsmath}
\usepackage{mathrsfs}
\usepackage[misc]{ifsym}
\usepackage{diagbox}
\usepackage{graphicx} 
\usepackage{amsmath,latexsym,amssymb,amsfonts,amsbsy}
\usepackage{algorithm}
\usepackage{algpseudocode}
\usepackage{color,xcolor}
\usepackage{booktabs}
\usepackage{threeparttable}
\usepackage{multicol}
\usepackage{multirow}
\usepackage{epstopdf}
\numberwithin{equation}{section}
\newtheorem{theorem}{Theorem}[section]
\newtheorem{lemma}[theorem]{Lemma}

\newtheorem{definition}[theorem]{Definition}

\newtheorem{corollary}[theorem]{Corollary}
\newtheorem{proposition}[theorem]{Proposition}
\newtheorem{condition}[theorem]{Condition}

\begin{document} 
\begin{sloppypar}
    \begin{frontmatter}
        \title{Numerical Analysis of Penalty--based Ensemble Methods}
        \author{Rui Fang\corref{cor1}}
\ead{ruf10@pitt.edu}
\ead[url]{https://ruf10.github.io}
\cortext[cor1]{Corresponding author: Rui Fang}

\begin{abstract}

The chaotic nature of fluid flow and the uncertainties in initial conditions limit predictability. Small errors that occur in the initial condition can grow exponentially until they saturate at $\mathcal{O}$(1). Ensemble forecasting averages multiple runs with slightly different initial conditions and other data to produce more accurate results and extend the predictability horizon. However, they can be computationally expensive. We develop a penalty--based ensemble method with a shared coefficient matrix to reduce required memory and computational cost and thereby allow larger ensemble sizes.
   Penalty methods relax the incompressibility condition to decouple the pressure and velocity, %Penalty methods relax incompressibility to uncouple the velocity and pressure,
   reducing memory requirements. This report gives stability proof and an error estimate of the penalty--based ensemble method, extends it to the Navier--Stokes equations with random variables using Monte Carlo sampling, and validates the method's accuracy and efficiency with three numerical experiments.
        \end{abstract}

    \end{frontmatter}

\section{Introduction}
\label{sec:1}
Unstable systems have finite predictability horizons, Lorenz \cite{lorenz1963deterministic, lorenz1963predictability}. The chaotic nature of fluid flow and the uncertainties in initial conditions limit predictability. Under different initial conditions, the trajectories of the flow spread. Small errors in the (uncertain) initial conditions can grow exponentially until $\mathcal{O}(1)$, resulting in a loss of prediction ability \cite{lorenz1985growth}.\\ 
Ensemble methods address the uncertainty in problem data by conducting numerical simulations with various initial and boundary conditions, external forces, and other data, % different initial and boundary conditions, body forces, and other data,
Kalnay \cite{kalnay2003atmospheric}. Monte Carlo forecasting with a sample size as small as $8$ will provide the best estimate, the ensemble mean, Leith \cite{leith1974theoretical}. Assume we have an ensemble of size \(J\). At each timestep, the execution process needs to assemble and solve \(J\) separate linear systems. Ensemble methods offer improved predictability but are computationally expensive. 

We develop a penalty--based ensemble method to reduce the computational cost and address the predictability limitations of flows. The method uses a shared coefficient matrix with different right--hand side vectors and relaxes the incompressibility condition to reduce the space complexity of the model while maintaining accuracy. Further, savings in memory and operations are obtained by eliminating the $J$ pressure variables. In this report, we derive a stability proof and an error estimate and conduct three numerical tests to validate the method. In addition, we extend it to the Navier--Stokes equations (NSE) with random variables. 

The incompressible NSE is given by
\begin{equation}\label{nse}
\begin{gathered}
\frac{\partial u}{\partial t} + u\cdot \nabla u - \nu \Delta u +\nabla p = f(x,t),\text{ and }
\nabla \cdot u =0,
\end{gathered}
\end{equation}
where $u$ denotes the flow velocity and $p$ denotes the flow pressure. The viscosity is denoted by $\nu$, and $f$ is the body force. In equation (\ref{nse}), the pressure is a Lagrange multiplier to enforce the incompressibility constraint, E and Liu \cite{weinan1995projection}. 
The penalty method relaxes incompressibility by 
replacing
\begin{equation*}
\nabla \cdot u=0 \text{ with } \nabla \cdot u^\epsilon+\epsilon p^\epsilon=0, \text{ for } 0< \epsilon << 1,
\end{equation*}
and hence it uncouples $u$ and $p$ and yields the penalized NSE:
\begin{equation}\label{penalty_conti}
\begin{gathered}
\frac{\partial u^\epsilon}{\partial t} + u^\epsilon \cdot \nabla u^\epsilon +\frac{1}{2} (\nabla \cdot u^\epsilon) u^\epsilon - \nu \Delta u^\epsilon +\nabla p^\epsilon = f,\\
\nabla \cdot u^{\epsilon}+\epsilon p^{\epsilon} =0, \text{ where } 0<\epsilon << 1.
\end{gathered}
\end{equation}
One can eliminate the pressure by setting $p^\epsilon=-\frac{1}{\epsilon}\nabla \cdot u^\epsilon$.

We adopt the ensemble approach of Nan and Layton \cite{jiang2014algorithm} to the penalized NSE, using a shared coefficient matrix with different right--hand sides. We suppress the spatial discretization to present the idea. We define the ensemble mean and fluctuation at the timestep $t_n$:
\begin{equation*}
\begin{gathered}
\langle u^\epsilon\rangle^n: = \frac{1}{J} \sum_{j=1}^J u^{\epsilon,n}_j, \text{ and } U^{\epsilon,n}_j:= u^{\epsilon,n}_j - \langle u^\epsilon\rangle^n,
\end{gathered}
\end{equation*}
where $u^{\epsilon,n}_{j}$ is the penalized velocity for the $j_{th}$ ensemble member. We use an implicit--explicit time discretization which allows the coefficient matrix to be independent of the ensemble member, which yields the following:
find $u^{\epsilon,n+1}_j$ in the velocity space and  $p^{\epsilon,n+1}_j$ in the pressure space: 
\begin{equation}\label{penalty_ensemble}
\begin{gathered}
\frac{u^{\epsilon, n+1}_j-u^{\epsilon, n}_j}{\Delta t} + \langle u^\epsilon\rangle^n\cdot \nabla u^{\epsilon,n+1}_j + \frac{1}{2} (\nabla \cdot \langle u^\epsilon\rangle^n) u^{\epsilon,n+1}_j \\
+U^{\epsilon, n}_j\cdot \nabla u^{\epsilon,n}_j + \frac{1}{2} (\nabla \cdot U^{\epsilon, n}_j) u^{\epsilon,n}_j
-\nu\Delta u_j^{\epsilon, n+1}
+\nabla p^{\epsilon, n+1}_{j} = f_j^{n+1},\\
\nabla\cdot u_j^{\epsilon, n+1} +\epsilon p_j^{\epsilon, n+1}= 0.
\end{gathered}
\end{equation}  
Here $\epsilon$ is the same for all ensemble members to ensure a shared coefficient matrix. The ensemble mean drives the flow. We can eliminate the pressure by setting $p^{\epsilon,n+1}_j=-\frac{1}{\epsilon} \nabla \cdot u^{\epsilon,n+1}_j$ to reduce the memory. 
%related work
\subsection{Related work}
Epstein \cite{epstein1969stochastic} introduced the first forecasting method that explicitly accounted for the uncertainty in atmospheric model predictions, known as the stochastic--dynamics forecasting method, in 1969. Leith \cite{leith1974theoretical} later proposed using ensemble forecasting with multiple members instead of a single realization. He showed that the ensemble mean from Monte Carlo ensembles can achieve accurate results without linear regression. 
Luo and Wang \cite{luo2018ensemble} studied an ensemble algorithm for the deterministic and random parabolic partial differential equations which led to a single discrete system with multiple right--hand side vectors. 

Temam \cite{temam1968methode} first introduced the penalty method with a modified nonlinear term to ensure energy dissipation. He proved in \cite{temam1968methode} that $\lim_{\epsilon \to 0} (u^{\epsilon}, p^{\epsilon})= (u,p)$.
Penalty methods have been widely studied, including Falk \cite{falk1976finite}, Shen \cite{shen1995error} and He \cite{he2010penalty}, He and Li \cite{he2005optimal}. We can speed up the calculation by eliminating the pressure by $p^\epsilon = -\frac{1}{\epsilon} \nabla \cdot u^\epsilon$, Heinrich and Vionnet \cite{heinrich1995penalty}. The error of velocity depends on the penalty parameter $\epsilon$, as shown by Bercovier and Engelman (1979) \cite{bercovier1979finite}. The condition number of the penalized system was studied in Layton and Xu \cite{layton2023conditioning}, Hughes and Liu and Brooks \cite{hughes1979finite}. Adapting penalty parameters, exploiting $\epsilon$--sensitivity, can help with ill--conditioning and provide better accuracy \cite{kean2023doubly, xie2022adaptive, fang2024numerical,layton2020doubly}, and pressure recovery in \cite{kean2020error}. Some preliminary tests of the penalty--based ensemble method are studied in Fang \cite{fang2023penalty}.

\section{Notations and preliminaries} 

Let $D\subset \mathbb{R}^d$ be an open regular domain, where $d=2,3$. The $L^2(D)$ norm is denoted as $\|\cdot\|$, and the inner product is denoted as $(\cdot,\cdot)$. Similarly, we define the $L^p(D)$ norms $\|\cdot\|_{L^p}$,  and the Sobolev $W^{k}_{p}(D)$ norms $\|\cdot\|_{W^{k}_{p}}$. We denote the Sobolev space $W^{k}_2(D)$ with norm $\|\cdot\|_{k}$ as $H^{k}(D)$. We define the norms for the functions $v(x,t)$ defined on $(0, T)$, for $1\leq m<\infty$,
\begin{equation}
\begin{gathered}
\|v\|_{\infty,k} := EssSup_{[0,T]} \|v(t,\cdot)\|_{k},\
\|v\|_{m,k} := \left(\int_{0}^T\|v(t,\cdot)\|^{m}_{k}\, dt \right)^{1/m}.
\end{gathered}
\end{equation}
The discrete--time equivalents of the norms are denoted as follows:
\begin{equation}
\begin{gathered}
\||v|\|_{\infty, k} := \max_{0 \leq n \leq N} \|v^n\|_k,  \text{ and }
\||v|\|_{m,k} := \left(\sum_{n=0}^{N}\|v^n\|^m_k \Delta t \right)^{1/m}.
\end{gathered}
\end{equation}
Let $(\Omega, \mathcal{F}, P)$ be a complete probability space, where $\Omega$ is the set of outcomes, $\mathcal{F} \subset 2^{\Omega}$ is the $\sigma$--algebra of events, and $P: \mathcal{F} \to [0,1]$ is a probability measure. We denote the set of all integrable functions for the probability measure $P$ as the $L^1_P(\Omega)$. Suppose a random variable $Y$ such that $Y
\in L^1_P(\Omega)$, we define the expected value of $Y$ as follows:
\begin{equation*}
E[Y] = \int_{\Omega} Y(\omega) \, d P(\omega).
\end{equation*}
The stochastic Sobolev spaces are denoted by 
\begin{equation*}\widetilde{W}^{k}_p:= L^p_P\left(\Omega, W^{k}_p(D)\right). 
\end{equation*}
$\widetilde{W}^{k}_p$ contains stochastic functions $v:\Omega \times D \to \mathbb{R}$, that are measurable with respect to the product $\sigma$--algebra $\mathcal{F}\otimes \mathcal{B}(D)$, where $\mathcal{B}$ is a Borel set. $\widetilde{W}^{k}_p$ is equipped with the averaging norms
\begin{equation*}
\|v\|_{\widetilde{W}^{k}_p} = \left(E\bigg[\|v\|^p_{{W}^{k}_p(D)}\bigg]\right)^{1/p}.
\end{equation*}
Note when $p=2$, the above space is a Hilbert space and we write $\widetilde{W}^{k}_2(D) = \widetilde{H}^{k}(D)$.

%We assume the following regularity:
%\begin{equation*}
%\begin{gathered}
%u_j \in L^\infty (0,T; H^{k+1}(D)) \cap H^1(0,T; H^{k+1}(D)) \cap H^2(0,T; L^2(D))\\
%p_j \in L^2(0,T; H^{s+1}(D)), \text{ and } f_j \in L^\infty(0,T; L^{2}(D)).
%\end{gathered} 
%\end{equation*}\\
%Also, additional regularity assumption from Shen \cite{shen1995error}.
%\begin{equation}\label{shen_assumption_1995}
%t f_t \in L^2(0,T; L^2(D))
%\end{equation}
%to prove
%\begin{equation*}
%tp_t \in L^2(0,T;H^1(D)).
%\end{equation*}

\begin{lemma} (\text{See Layton \cite{layton2008introduction}, p. 28, p. 29})\label{Poincare_Friedriches_inequality}
Suppose $\Gamma_0 \subset \partial D$ has a positive measure. Let
 \begin{equation}
H^{1}_0(D): = \{v \in L^2(D): \nabla v\in L^2(D) \text{ and } v =0 \text{ on } \Gamma_0\}.
 \end{equation}
Then, there is a positive constant $C_{PF}$ such that
 \begin{equation}
 \|v\|\leq C_{PF}\|\nabla v\| \text{ for every } v\in H^{1}_0(D).
 \end{equation}
 Thus, $\|\nabla v\|$ and $\|v\|$ are equivalent norms on $H^{1}_0(D)$.
\end{lemma}
The space $H^{-k}(D)$ is the dual space of bounded linear functionals on $H^{k}_0(D)$.
A norm for $H^{-1}(D)$ is given by
\begin{equation}
\|f\|_{-1} = \sup_{0 \neq v \in H^{1}_{0}(D)} \frac{(f,v)}{\|\nabla v\|}.
\end{equation}

%\begin{lemma} For any $f$,
%$\|f\|_{-1} \leq \|f\|$.
%\end{lemma}
%\begin{proof}
%\begin{equation}
%\|f\|_{-1} = \sup_{0 \neq v \in H^{1}_{0}(D)} \frac{(f,v)}{\|\nabla v\|}\leq \sup_{0 \neq v \in H^{1}_{0}(D)} \frac{\|f\|\|v\|}{\|\nabla v\|}.
%\end{equation}
%By Lemma \ref{Poincare_Friedriches_inequality},
%\begin{equation*}
%\|f\|_{-1}\leq \sup_{0 \neq v \in H^{1}_{0}(D)} \frac{C_{PF}\|f\|\|\nabla v\|}{\|\nabla v\|}=C_{PF}\|f\|.
%\end{equation*}
%\end{proof}
\begin{lemma}(See Layton \cite{layton2008introduction}, p. 11 )
Let $D\subset \mathbb{R}^2$ or $\mathbb{R}^3$. If $f\in L^2({D})$, then
\begin{equation*}
\|f\|_{-1}\leq C_{PF} \|f\|< \infty.
\end{equation*}
\end{lemma}
Let $X$ be the velocity space and $Q$ be the pressure space:
\begin{equation}
\begin{gathered}
X:= (H^{1}_{0}(D))^d, \text{ and }Q:= L^2_0(D).
\end{gathered}
\end{equation}
We denote the conforming velocity and pressure finite element spaces as follows:
\begin{equation*}
X^h \subset X \text{ and } Q^h \subset Q.
\end{equation*} We assume that $(X^h,Q^h )$ satisfies the following approximation properties and the Ladyzhenskaya-Babushka-Brezzi Conditon ($LBB^h$). For $u\in H^{m+1}(D)^d$ and $p\in H^m(D)$,
\begin{equation}\label{interpolation_inequality}
\begin{gathered}
\inf_{v\in X^h} \|\nabla (u-v)\|\leq C h^m |u|_{m+1},\\ 
\inf_{v\in X^h} \|u-v\| \leq C h^{m+1}|u|_{m+1},\\
\inf_{q\in Q^h}\|p-q\| \leq C h^m |p|_m.
\end{gathered}
\end{equation}
\begin{condition}(See Layton \cite{layton2008introduction} p. 62, $LBB^h$ condition)\label{LBB_h}
Suppose $(X^h, Q^h)$ satisfies
\begin{equation}
\inf_{ q^h \in Q^h} \sup_{v_h \in X^h} \frac{(q^h,\nabla \cdot v_h)}{\|v_h\|\|q^h\|}\geq \beta^h >0,
\end{equation}
where $\beta^h$ is bounded away from zero uniformly in $h$.
\end{condition}
Condition \ref{LBB_h} is equivalent to
\begin{equation*}
\beta^h \|q^h\|\leq \sup_{v_h \in X^h} \frac{(q^h, \nabla \cdot v_h)}{\|v_h\|}.
\end{equation*}

We assume the mesh with quasi--uniform triangulation and finite element spaces satisfy the inverse inequality: 
\begin{equation}
\begin{gathered}
h\|\nabla v_h\|\leq C \|v_h\|  \ \ \forall  v_h \in X^h.
\end{gathered}
\end{equation}
\begin{lemma}(See Ladyshenskaya \cite{ladyzhenskaya1969mathematical})\label{Ladyzhenskaya} For any vector function $u: \mathbb{R}^d \to \mathbb{R}^d$ with compact support and with finite $L^p$ norms:
\begin{equation}
\begin{gathered}
\quad \|u\|_{L^4(\mathbb{R}^2)} \leq 2^{1/4} \|u\|^{1/2}_{L^2(\mathbb{R}^2)}\|\nabla u\|^{1/2}_{L^2(\mathbb{R}^2)}, (d=2),\\
\quad \|u\|_{L^4(\mathbb{R}^3)} \leq \frac{4}{3\sqrt{3}}\|u\|^{1/4} \|\nabla u\|^{3/4}, (d=3),\\
\quad \|u\|_{L^6(\mathbb{R}^3)}\leq \frac{2}{\sqrt{3}}\|\nabla u\|, (d =3).
\end{gathered}
\end{equation}
\end{lemma}
\begin{lemma}(A discrete Gronwall lemma, see Lemma 5.1, p. 369, \cite{heywood1990finite} ) \label{discrete_gronwall} Let $\Delta t, B, a_n, b_n, c_n ,d_n$ be
nonnegative numbers such that for $l
\geq 1$:
\begin{equation}
\begin{gathered}
a_l + \Delta t \sum_{n=0}^l b_n\leq \Delta t \sum_{n=0}^{l-1} d_n a_n + \Delta t \sum_{n=0}^{l}c_n + B, \text{ for } l\geq 0,
\end{gathered}
\end{equation}
then for all $\Delta t>0$,
\begin{equation}
\begin{gathered}
a_l + \Delta t \sum_{n=0}^l b_n\leq \exp(\Delta t \sum_{n=0}^{l-1} d_n) (\Delta t \sum_{n=0}^l c_n +B).
\end{gathered}
\end{equation}
\end{lemma}

\begin{lemma}\label{holderandyoung} (\text{See Layton \cite{layton2008introduction}, p. 7,}
H\"older's and Young's inequalities)
For any $\xi>0$, $1\leq p<\infty$, and $\frac{1}{p}+\frac{1}{q}=1$, the H\"older and Young's inequalities:
\begin{equation*}
\begin{gathered}
    (u,v)\leq \|u\|_{L^p}\|v\|_{L^q},\
    (u,v)\leq \frac{\xi}{p} \|u\|_{L^p}^p+
    \frac{\xi^{-q/p}}{q} \|v\|_{L^q}^q.
    \end{gathered}
\end{equation*}
The generalization with three functions,
\begin{equation}
\begin{gathered}
|f g h|\leq \|f\|_{L^p}\|g\|_{L^q} \|h\|_{L^r}, \text{ where } \frac{1}{p}+ \frac{1}{q}+ \frac{1}{r} =1.
\end{gathered}
\end{equation}
\end{lemma}
The standard skew--symmetric trilinear form is $\forall u,v,w \in X$,
\begin{equation*}
\begin{split}
b^*(u,v,w) :=\frac{1}{2}(u\cdot \nabla v, w)-\frac{1}{2} (u\cdot \nabla w, v). 
\end{split}
\end{equation*}

\begin{lemma}\label{trilinear_upper_bound} (See Layton \cite{layton2008introduction}, p. 11, Lemma 3) For any $u, v, w\in X$, there is $C=C(D)$ such that
\begin{equation*}
\begin{gathered}
\left| \int_{D} u \cdot \nabla v \cdot w\, dx \right|\leq C \|\nabla u\| \|\nabla v\| \|\nabla w\|, \text{ and }\\
\left|\int_{D} u\cdot \nabla v \cdot w \, dx \right|\leq C\|u\|^{\frac{1}{2}}\|\nabla u\|^{\frac{1}{2}}\|\nabla v\| \|\nabla w\|.
\end{gathered}
\end{equation*}
\end{lemma}
\begin{lemma} \label{another_form_of_trilinear}(\text{See Layton \cite{layton2008introduction}, p. 123, p. 155}) $\forall u,v,w \in X$,
\begin{equation*}
\begin{split}
b^*(u,v,w)=(u\cdot \nabla v, w)+\frac{1}{2} \left((\nabla \cdot u) v, w\right).
\end{split}
\end{equation*}
\end{lemma}
\begin{lemma} (See Layton \cite{layton2008introduction} and Girault and Raviart \cite{girault2012finite}) $b^*(u,v,w)$ satisfies the following bounds:
%\begin{equation}
%\begin{gathered}
%b^*(u,v,w)\leq C \sqrt{\|u\|\|\nabla u\|}\|\nabla v\| \|\nabla w\|,\\
%b^*(u,v,w)\leq C \|\nabla u\|\|\nabla v\| \sqrt{\|w\|\|\nabla w\|},\\
%b^*(u,v,w)\leq C \|\nabla u\|\|\nabla v\|\|\nabla w\|.
%\end{gathered}
%\end{equation}
\begin{equation}
b^*(u,v,w) \leq \begin{cases} 
C \sqrt{\|u\|\|\nabla u\|}\|\nabla v\| \|\nabla w\|, & \\
C \|\nabla u\|\|\nabla v\| \sqrt{\|w\|\|\nabla w\|}, & \\
C \|\nabla u\|\|\nabla v\|\|\nabla w\|. 
\end{cases}
\end{equation}
for all $u,v,w \in X$.
\end{lemma}
\begin{definition} $P_{Q^h}$ is the $L^2$ projection of $Q$ onto $Q^h$. That is, for any $q\in Q$, $P_{Q^h}(q)$ satisfies
\begin{equation*}
(P_{Q^h} (q)-q , q^h)= 0,
\ \forall q^h \in Q^h.
\end{equation*}
\end{definition}
\section{Penalty--based ensemble method}
We define the final time $T$ and timstep size at the $n_{th}$ step $\Delta t_n$. The total number of steps $N$ is given by $N = T / \Delta t$. The fully--discrete approximation is then given $(u^{\epsilon,n}_{j,h},p^{\epsilon,n}_{j,h})  \in (X^h, Q^h)$, find $(u^{\epsilon,n+1}_{j,h}, p^{\epsilon, n+1}_{j,h} ) \in (X^h, Q^h)$ satisfying: 
\begin{equation}\label{weak_form_method}
\begin{gathered}
\frac{1}{\Delta t_n}(u^{\epsilon, n+1}_{j,h}-u^{\epsilon,n}_{j,h}, v_h) + b^*(\langle u^\epsilon_h\rangle^n, u^{\epsilon,n+1}_{j,h},v_h)+ b^*(u^{\epsilon,n}_{j,h} -\langle u^\epsilon_h\rangle^n, u^{\epsilon,n}_{j,h},v_h)\\
+\nu (\nabla u^{\epsilon,n+1}_{j,h}, \nabla v_h) -(p^{\epsilon,n+1}_{j,h}, \nabla \cdot v_h) +(q^h, \nabla \cdot u^{\epsilon, n+1}_{j,h}) + \epsilon (p^{\epsilon, n+1}_{j,h}, q^h)= (f_j^{n+1},v_h),
\end{gathered}
\end{equation}
for all $(v_h, q^h)\in (X^h,Q^h)$. 

Due to the stretching term $b^*(u^{\epsilon,n}_{j,h} -\langle u^\epsilon_h\rangle^n, u^{\epsilon,n}_{j,h},v_h)$, we need the CFL timestep restriction:
\begin{equation}\label{cfl}
C\frac{\Delta t}{\nu h} \|\nabla (u^{\epsilon,n}_{j,h}- \langle u^\epsilon_{h}\rangle^n)\|^2 \leq 1. 
\end{equation}
If equation (\ref{cfl}) is satisfied, we proceed to the next timestep. Otherwise, we halve the timestep and repeat the current step.
%
%\begin{algorithm}
%\caption{}\label{alg:cap}
%\begin{algorithmic}
%\Require $\Delta t$, $\epsilon$, $u^n$, and $U^{\epsilon, n}_j, j=1,\ldots, J$ at time $t$.
%\While{$t <  T$}:
%\State Compute the average $\langle u \rangle^n = \frac{1}{J} \sum_{j=1}^J U^{\epsilon, n}_j$.
%\State Solve the next step velocity $u^{\epsilon,n+1}_j$ by the method in equation (\ref{weak_form_method}).
%\State Update the average $\langle u \rangle^{n+1} = \frac{1}{J} \sum_{j=1}^J u^{\epsilon,n+1}_j$.
%\State Calculate $\max_{j} \|\nabla (u^{\epsilon,n+1}_j - \langle u \rangle^{n+1})\|$ and verify the CFL condition.
%\If{the CFL condition is satisfied}
%    \State update time by $t = t + \Delta t$.
%  \Comment{Succeed, and go to next timestep}
%\Else
%    \State $\Delta t = \frac{1}{2} \Delta t$.
%\EndIf
%\EndWhile
%\end{algorithmic}
%\end{algorithm}

%
\subsection{Stability }\label{sec: stability} 
Let the difference between the ensemble member $j$ and the ensemble average be denoted as
\begin{equation}
U^{\epsilon,n}_j := u^{\epsilon,n}_{j,h}-\langle u^{\epsilon}_h\rangle^n.
\end{equation}
In Theorem \ref{sec: stability}, we prove the method is nonlinearly and long--time energy stable under the CFL condition:
\begin{equation*}
C\frac{\Delta t}{\nu h} \|\nabla U^{\epsilon,n}_j\|^2 \leq 1.
\end{equation*}
Define
\begin{equation*}
p^{\epsilon}_{j,h} =-\frac{1}{\epsilon}  P_{Q^h}(\nabla \cdot u^{\epsilon}_{j,h}).
\end{equation*}
\begin{theorem}
\label{Stability of BEFE-Ensemble penalty} Suppose the following timestep condition holds:
\begin{equation}\label{stability_condition}
C\frac{\Delta t}{\nu h} \|\nabla U^{\epsilon,n}_j\|^2\leq 1, j=1,\ldots, J.
\end{equation}
It yields that for any $N\geq 1$:
\begin{equation}\label{stability}
\begin{gathered}
\frac{1}{2} \|u^{\epsilon,N}_{j,h}\|^2 + \frac{1}{4}\sum_{n=0}^{N-1}\|u^{\epsilon,n+1}_{j,h}- u^{\epsilon,n}_{j,h}\|^2 + \frac{\nu \Delta t}{4} \|\nabla u^{\epsilon,N}_{j,h}\|^2 \\
+ \frac{\Delta t}{\epsilon}\sum_{n=0}^{N-1} \|P_{Q^h}(\nabla \cdot u^{\epsilon,n+1}_{j,h})\|^2 
+\frac{\nu \Delta t}{4}\sum_{n=0}^{N-1} \|\nabla u^{\epsilon,n+1}_{j,h}\|^2 \\
\leq \frac{\Delta t}{2 \nu}\sum_{n=0}^{N-1} \|f^{n+1}_{j,h}\|^2_{-1} + \frac{1}{2} \|u^{\epsilon,0}_{j,h}\|^2 +\frac{\nu \Delta t}{4} \|\nabla u^{\epsilon,0}_{j,h}\|^2.
\end{gathered}
\end{equation}
\end{theorem}
\begin{proof}
We write $p^{\epsilon}_{j,h} =-\frac{1}{\epsilon}  P_{Q^h}(\nabla \cdot u^{\epsilon}_{j,h})$ in the momentum equation, and inner product with $v_h\in X^h$. It yields
\begin{equation}\label{weak-form}
\begin{gathered}
\frac{1}{\Delta t_n}(u^{\epsilon, n+1}_{j,h}-u^{\epsilon,n}_{j,h}, v_h) + \nu (\nabla u^{\epsilon,n+1}_{j,h}, \nabla v_h) +\frac{1}{\epsilon}(P_{Q^h}(\nabla \cdot u^{\epsilon,n+1}_{j,h}),\nabla \cdot v_h)\\
+b^*(\langle u^\epsilon_h\rangle^n, u^{\epsilon,n+1}_{j,h},v_h)+ b^*(u^{\epsilon,n}_{j,h} -\langle u^\epsilon_h\rangle^n, u^{\epsilon,n}_{j,h},v_h)
 = (f_j^{n+1},v_h).
\end{gathered}
\end{equation}
Set $v_h = u^{\epsilon, n+1}_{j,h}$. Multiply $\Delta t$ to both sides of the equation (\ref{weak-form}) and apply the polarization identity:
\begin{equation}
\begin{gathered}
\frac{1}{2}\|u^{\epsilon, n+1}_{j,h}\|^2 -\frac{1}{2}\|u^{\epsilon,n}_{j,h}\|^2 +\frac{1}{2} \|u^{\epsilon,n+1}_{j,h}- u^{\epsilon,n}_{j,h}\|^2 + \Delta t b^*(U^{\epsilon,n}_j,u^{\epsilon,n}_{j,h}, u^{\epsilon, n+1}_{j,h})\\
+\nu \Delta t \|\nabla u^{\epsilon,n+1}_{j,h}\|^2 +  \frac{\Delta t}{\epsilon} \|P_{Q^h}(\nabla \cdot u^{\epsilon,n+1}_{j,h})\|^2 = \Delta t (f^{n+1}_j, u^{\epsilon,n+1}_{j,h}).
\end{gathered}
\end{equation}
Apply Young's inequality to $(f^{n+1}_j, u^{\epsilon,n+1}_{j,h})$ gives:
\begin{equation}
    \begin{gathered}
\frac{1}{2}\|u^{\epsilon,n+1}_{j,h}\|^2 -\frac{1}{2}\|u^{\epsilon,n}_{j,h}\|^2 +\frac{1}{2} \|u^{\epsilon,n+1}_{j,h}- u^{\epsilon,n}_{j,h}\|^2 + \Delta t b^*(U^{\epsilon,n}_j,u^{\epsilon,n}_{j,h},u^{\epsilon,n+1}_{j,h})\\
+\nu \Delta t \|\nabla u^{\epsilon,n+1}_{j,h}\|^2 +  \frac{\Delta t}{\epsilon} \|P_{Q^h}(\nabla \cdot u^{\epsilon,n+1}_{j,h})\|^2 \leq \frac{\nu \Delta t}{2}\| \nabla u^{\epsilon,n+1}_{j,h}\|^2+\frac{\Delta t}{2\nu} \|f_j^{n+1}\|^2_{-1}.
    \end{gathered}
\end{equation}
Next, we treat the trilinear term with the help of inverse inequalities and interpolation, 
\begin{equation}
\begin{gathered}
-\Delta t b^*(U^{\epsilon,n}_j,u^{\epsilon,n}_{j,h}, u^{\epsilon,n+1}_{j,h})
=-\Delta t b^*(U^{\epsilon,n}_j,u^{\epsilon,n}_{j,h}, u^{\epsilon,n+1}_{j,h}-u^{\epsilon,n}_{j,h})\\
\leq C\Delta t \|\nabla U^{\epsilon,n}_j\| \|\nabla u^{\epsilon,n}_{j,h}\| \left(\|\nabla (u^{\epsilon,n+1}_{j,h}- u^{\epsilon,n}_{j,h})\| \|u^{\epsilon,n+1}_{j,h}- u^{\epsilon,n}_{j,h}\|\right)^{1/2}\\
\leq C\Delta t \|\nabla U^{\epsilon,n}_j\| \|\nabla u^{\epsilon,n}_{j,h}\| \frac{1}{\sqrt{h}} \|u^{\epsilon,n+1}_{j,h}- u^{\epsilon,n}_{j,h}\| \\
\leq C \frac{\Delta t^2}{h} \|\nabla U^{\epsilon,n}_j\|^2 \|\nabla u^{\epsilon,n}_{j,h}\|^2 + \frac{1}{4} \|u^{\epsilon,n+1}_{j,h}- u^{\epsilon,n}_{j,h}\|^2.
\end{gathered}
\end{equation}
Combine terms, we have
\begin{equation}
\begin{gathered}
\frac{1}{2}\|u^{\epsilon,n+1}_{j,h}\|^2 -\frac{1}{2}\|u^{\epsilon,n}_{j,h}\|^2 +\frac{1}{4} \|u^{\epsilon,n+1}_{j,h}- u^{\epsilon,n}_{j,h}\|^2
+\frac{\nu \Delta t}{2}\|\nabla u^{\epsilon,n+1}_{j,h}\|^2 \\
+  \frac{\Delta t}{\epsilon} \| P_{Q^h} (\nabla \cdot u^{\epsilon,n+1}_{j,h})\|^2
\leq \frac{\Delta t}{2\nu} \|f_j\|^2_{-1} + C \frac{\Delta t^2}{h} \|\nabla U^{\epsilon,n}_j\|^2 \|\nabla u^{\epsilon,n}_{j,h}\|^2.
\end{gathered}
\end{equation}
Add and subtract $\frac{\nu \Delta t}{4} \|\nabla u^{\epsilon,n}_{j,h}\|^2$, we have
\begin{equation}\label{n_step}
 \begin{gathered}
 \frac{1}{2}\|u^{\epsilon,n+1}_{j,h}\|^2 -\frac{1}{2}\|u^{\epsilon,n}_{j,h}\|^2 +\frac{1}{4} \|u^{\epsilon,n+1}_{j,h}- u^{\epsilon,n}_{j,h}\|^2+\frac{
 \nu\Delta t}{4}\|\nabla u^{\epsilon,n+1}_{j,h}\|^2\\
 +\frac{\nu \Delta t}{4}\left(\|\nabla u^{\epsilon,n+1}_{j,h}\|^2 -\|\nabla u^{\epsilon,n}_{j,h}\|^2\right)
+\frac{\Delta t}{\epsilon} \|P_{Q^h}(\nabla \cdot u^{\epsilon,n+1}_{j,h})\|^2\\
 + \frac{\nu\Delta t}{4}(1- \frac{C \Delta t}{h} \|\nabla U^{\epsilon,n}_j\|^2) \|\nabla u^{\epsilon,n}_{j,h}\|^2
 \leq \frac{\Delta t}{2 \nu} \|f^{n+1}_j\|^2_{-1}.
 \end{gathered}
\end{equation}
With the CFL condition in equation (\ref{stability_condition}), equation (\ref{n_step}) reduces to:
\begin{equation}
\begin{gathered}
\frac{1}{2}\|u^{\epsilon,n+1}_{j,h}\|^2 -\frac{1}{2}\|u^{\epsilon,n}_{j,h}\|^2 +\frac{1}{4} \|u^{\epsilon,n+1}_{j,h}- u^{\epsilon,n}_{j,h}\|^2+\frac{\nu \Delta t}{4}(\|\nabla u^{\epsilon, n+1}_{j,h}\|^2 -\|\nabla u^{\epsilon, n}_{j,h}\|^2)\\
 + \frac{\Delta t}{\epsilon} \|P_{Q^h}(\nabla \cdot u^{\epsilon,n+1}_{j,h})\|^2 
 + \frac{\nu\Delta t}{4} \|\nabla u^{\epsilon,n+1}_{j,h}\|^2  \leq \frac{\Delta t}{2 \nu} \|f^{n+1}_j\|^2_{-1}.
\end{gathered}
\end{equation}
Sum over all 
$n$ from $0$ to 
$N-1$, we have the final result.
\end{proof}
\begin{lemma} (See Evans and Rosenthal \cite{evans2004probability}, p. 149, Theorem 3.3.1 ) \label{ensemble_equality} Let $u_j$ be the $j_{th}$ ensemble member, and $\langle u \rangle=\frac{1}{J}\sum_{j=1}^J u_j$. Then the variance is equal to the
second moment minus the square of the first moment.
\begin{equation}
\frac{1}{J}\sum_{j=0}^J\|u_j-\langle u \rangle\|^2 = \|\langle u \rangle\|^2- \frac{1}{J}\sum_{j=0}^J \|u_j\|^2 .
\end{equation}
\end{lemma}
\begin{proposition}\label{average_equality}
\begin{equation}
\Delta t\sum_{n=0}^{N} \frac{1}{J} \sum_{j=1}^J  \|\nabla U^{\epsilon, n}_j\|^2 <C.
\end{equation}
\end{proposition}
\begin{proof}
By Lemma \ref{ensemble_equality},  
\begin{equation*}
\frac{1}{J}\sum_{j=0}^J\|\nabla u^{\epsilon,n}_{j,h}\|^2 = \|\nabla \langle u^{\epsilon}_h\rangle^n\|^2 +\frac{1}{J}\sum_{j=0}^J\|\nabla U^{\epsilon,n}_j\|^2.
\end{equation*}
Sum from $n=0$ to $n=N$, and multiply by $\Delta t$:
\begin{equation}
\Delta t \sum_{n=0}^N \frac{1}{J}\sum_{j=1}^J\|\nabla u^{\epsilon,n}_{j,h}\|^2 = \Delta t \sum_{n=0}^N \|\nabla \langle u^{\epsilon}_h\rangle^n\|^2 + \Delta t \sum_{n=0}^N \sum_{j=1}^J  \|\nabla U^{\epsilon,n}_j\|^2.
\end{equation}
Since
$\Delta t \sum_{n=0}^N \|\nabla \langle u^{\epsilon}_h\rangle^n\|^2 \geq 0$, and $\Delta t \sum_{j=1}^J \sum_{n=0}^N \|\nabla U^{\epsilon,n}_j\|^2\geq 0$,
it is sufficient to show that $\Delta t \sum_{n=0}^N \frac{1}{J}\sum_{j=1}^J\|\nabla u^{\epsilon,n}_{j,h}\|^2$ is bounded by a finite number. By Theorem \ref{Stability of BEFE-Ensemble penalty}, for $j =1,\ldots, J$, we have
\begin{equation*}
\Delta t \sum_{n=0}^N \|\nabla u^{\epsilon,n}_{j,h}\|^2< C.
\end{equation*}
Hence, we have
\begin{equation*}
\Delta t \sum_{n=0}^N \frac{1}{J}\sum_{j=1}^J\|\nabla u^{\epsilon,n}_{j,h}\|^2<\infty.
\end{equation*}
\end{proof}
\subsection{Error estimates}\label{ sec: error-analysis}
\begin{definition}\label{stokes_projection} Define
$P_s:(X,Q) \to (X^h, Q^h)$, the Stokes projection. $P_s (u,p)= (\Tilde{u}, \Tilde{p})$ satisfies: $\forall v_h \in X^h$ and $q^h \in Q^h$,
\begin{equation}
\begin{gathered}
    \nu(\nabla (u-\Tilde{u}),\nabla v_h)-(p-\Tilde{p}, \nabla \cdot v_h)=0, \\
    (\nabla \cdot (u-\Tilde{u}), q^h)=0.
\end{gathered}
\end{equation}
\end{definition}
\begin{proposition}\label{stokes-estimate} (See John \cite{john2016finite}, p. 164, Lemma 4.43)
Let the domain $D$ be bounded with polyhedral and Lipschitz continuous boundary and $(u,p)\in (X, Q)$. Suppose $LBB^h$ Condition \ref{LBB_h} holds, then it yields
\begin{equation}
\begin{gathered}
\|\nabla (u-\Tilde{u})\| \leq 2\left( 1+\frac{1}{\beta^h}\right) \inf_{v_h\in X^h} \|\nabla (u-v_h)\|+ \inf_{q^h \in Q^h}\|p-q^h\|,\\
\|p-\Tilde{p}\|\leq \frac{2}{\beta^h}\bigg\{ \left(1+\frac{1}{\beta^h}\right)\inf_{v_h\in X^h}\|\nabla (u-v_h)\|
+\inf_{ q^h \in Q^h} \|p-q^h\|  \bigg\}.\\
\end{gathered}
\end{equation}
\end{proposition}
Denote the error of the $j_{th}$ simulation at time $t_n$, $e^{\epsilon,n}_{j}:= u^{\epsilon,n}_j - u^{\epsilon, n}_{j,h}$. Here, $u^{\epsilon,n}_j$ is the solution of the penalized NSE at time $t_n$ and $u^{\epsilon,n}_{j,h}$ is the fully discretized solution of penalty--based ensemble method.

\begin{theorem}\label{convergence_ensemble_penalty}
Consider the method in equation (\ref{penalty_ensemble}) and assume the condition in equation (\ref{stability_condition}) holds for all $n$:
\begin{equation}
C\frac{\Delta t}{\nu h} \|\nabla U^{\epsilon,n}_j\|^2\leq 1, j=1,\ldots, J,
\end{equation}
then there are positive constants $C$ and $C_0$ independent of $h$ and $\Delta t$ such that:
\begin{equation}\label{convergence_eqn}
\begin{gathered}
\|e^{\epsilon, N}_{j,h}\|^2 + \frac{1}{2}\sum_{n=0}^{N-1} \|e^{\epsilon,n+1}_{j,h}- e^{\epsilon,n}_{j,h}\|^2 + \Delta t \nu \|\nabla e^{\epsilon,N}_{j,h}\|^2\\
+ C_0\Delta t \sum_{n=0}^{N-1}\nu\|\nabla e^{\epsilon,n+1}_{j,h}\|^2 \leq \exp (\alpha )\Biggl\{\|e^{\epsilon, 0}_{j,h}\|^2+ \Delta t \nu \|\nabla e^{\epsilon,0}_{j,h}\|^2\\
+ h^{2m}C(\nu) T\left( \||u^{\epsilon}_{j,t}|\|^2_{\infty,0} + \frac{1}{\nu^2} \||p^{\epsilon}_{j,t}|\|^2_{\infty,0} \right)+ (\Delta t)^3C(\nu) \||u^{\epsilon}_{j,tt}|\|^2_{\infty, 0}\\
+ h^{2m}\epsilon \Delta t C(\nu, \beta^h)\left(\||u^{\epsilon}_{j,t}|\|^2_{2,0}  + \||p^{\epsilon}_{j,t}|\|^2_{2,0} \right)\\
+ h^{2m}C(\nu) T\left(  \|| u^{\epsilon}_j|\|^2_{2,0} +\frac{1}{\nu^2} \||p^{\epsilon}_j|\|^2_{2,0}\right)+ C(\nu) (\Delta t)^2 \||\nabla u^{\epsilon}_{j,t}|\|^2_{\infty,0}
\Biggl\},
\end{gathered}
\end{equation}
where
\begin{equation*}
\alpha = C(\nu) \Delta t \sum_{n=0}^{N-1} \|\nabla u^{\epsilon, n+1}_j\|^4. 
\end{equation*}
\end{theorem}
\begin{proof}
We evaluate the continuous penalty--based NSE (equation (\ref{penalty_conti})) at time $t = t_{n+1}$. For any $v_h\in X^h$, and $q^h \in Q^h$, 
\begin{equation} \label{exact_weak_form}
    \begin{gathered}(\frac{u^{\epsilon,n+1}_{j}- u^{\epsilon,n}_j}{\Delta t} , v_h) + b^*(u^{\epsilon,n+1}_j, u^{\epsilon,n+1}_j, v_h)+\nu (\nabla u^{\epsilon,n+1}_j, \nabla v_h) \\
-(p^{\epsilon,n+1}_{j}, \nabla \cdot v_h) + (\nabla \cdot u^{\epsilon,n+1}_{j}, q^h)+ \epsilon(p^{\epsilon,n+1}_{j}, q^h)= (f_j^{n+1},v_h)-(r^{\epsilon,n+1}_j,v_h),
    \end{gathered}
\end{equation}
where 
\begin{equation*}
r^{\epsilon,n+1}_j = u^{\epsilon,n+1}_{j,t} -\frac{u^{\epsilon,n+1}_j - u^{\epsilon,n}_j}{\Delta t}.
\end{equation*}
Subtract equation (\ref{weak_form_method}) from equation (\ref{exact_weak_form}). We have
\begin{equation}\label{error_1}
\begin{gathered}
\frac{1}{\Delta t} (e^{\epsilon,n+1}_j-e^{\epsilon,n}_j,v_h) + b^*(u^{\epsilon,n+1}_j, u^{\epsilon,n+1}_j, v_h) - b^*(\langle u^{\epsilon}_h\rangle^n, u^{\epsilon,n+1}_{j,h}, v_h)\\
- b^*(u^{\epsilon,n}_{j,h}-\langle u^{\epsilon}_h\rangle^n, u^{\epsilon,n}_{j,h}, v_h) 
+\nu (\nabla e^{\epsilon,n+1}_j,\nabla v_h)-(p^{\epsilon,n+1}_{j} - p^{\epsilon,n+1}_{j,h}, \nabla \cdot v_h)\\
+(\nabla \cdot e^{\epsilon,n+1}_{j}, q^h)+\epsilon(p^{\epsilon,n+1}_{j}-p^{\epsilon,n+1}_{j,h},q^h) +(r^{\epsilon,n+1}_j, v_h)= 0.
\end{gathered}
\end{equation}
Let $\Tilde{u} \in X^h$ and $\Tilde{q}\in Q^h$, define $e^{\epsilon,n}_j = \eta^{\epsilon,n}_j -\phi^{\epsilon,n}_{j,h}$, where $\eta^{\epsilon,n}_j:=u^{\epsilon,n}_j - \Tilde{u}$, $\phi^{\epsilon,n}_{j,h}:= u^{\epsilon,n}_{j,h}-\Tilde{u}$.  
\begin{equation}\label{before-stokes} 
\begin{gathered}
\frac{1}{\Delta t} (\phi^{\epsilon,n+1}_{j,h}-\phi^{\epsilon,n}_{j,h}, v_h) +\nu (\nabla \phi^{\epsilon,n+1}_{j,h},\nabla v_h) -(p^{\epsilon,n+1}_{j,h}-\Tilde{q}, \nabla \cdot v_h )+(\nabla \cdot \phi^{\epsilon,n+1}_{j,h}, q^h)\\
+\epsilon(p^{\epsilon,n+1}_{j,h}-\Tilde{q}, q^h)
=\frac{1}{\Delta t}(\eta^{\epsilon, n+1}_j-\eta^{\epsilon, n}_j,v_h)+ \nu(\nabla \eta^{\epsilon, n+1}_j, \nabla v_h) -(p^{\epsilon,n+1}_{j}-\Tilde{q}, \nabla \cdot v_h ) \\
+(\nabla \cdot \eta^{\epsilon,n+1}_{j}, q^h)+\epsilon(p^{\epsilon,n+1}_{j}-\Tilde{q}, q^h)
+(r^{\epsilon,n+1}_j,v_h)\\
+ b^*(u^{\epsilon, n+1}_j, u^{\epsilon,n+1}_j, v_h)
-b^*(\langle u^{\epsilon}_h\rangle^n, u^{\epsilon,n+1}_{j,h}, v_h)
- b^*(U^{\epsilon,n}_j, u^{\epsilon,n}_{j,h}, v_h).
\end{gathered}
\end{equation}
Let $\Tilde{u} \in X^h$ and $\Tilde{q}\in Q^h$ satisfy the Stokes projection:
\begin{equation*}
\begin{gathered}
\nu( \nabla (u^{\epsilon,n+1}_{j}- \Tilde{u}),\nabla v_h) -(p^{\epsilon,n+1}_{j}-\Tilde{q}, \nabla \cdot v_h)=0 \text{ for all } v_h \in X^h,\\
 (\nabla \cdot (u^{\epsilon,n+1}_{j}- \Tilde{u}), q^h) =0\text{ for all } q^h \in Q^h.
\end{gathered}
\end{equation*}
Equation (\ref{before-stokes}) is simplified to
\begin{equation} 
\begin{gathered}
\frac{1}{\Delta t} (\phi^{\epsilon,n+1}_{j,h}-\phi^{\epsilon,n}_{j,h}, v_h) +\nu (\nabla \phi^{\epsilon,n+1}_{j,h},\nabla v_h) -(p^{\epsilon,n+1}_{j,h}-\Tilde{q}, \nabla \cdot v_h )+(\nabla \cdot \phi^{\epsilon,n+1}_{j,h}, q^h)\\
+\epsilon(p^{\epsilon,n+1}_{j,h}-\Tilde{q}, q^h)
=\frac{1}{\Delta t}(\eta^{\epsilon, n+1}_j-\eta^{\epsilon, n}_j,v_h)+\epsilon(p^{\epsilon,n+1}_{j}-\Tilde{q}, q^h)
+(r^{\epsilon,n+1}_j,v_h)\\
+ b^*(u^{\epsilon, n+1}_j, u^{\epsilon,n+1}_j, v_h)
-b^*(\langle u^{\epsilon}_h\rangle^n, u^{\epsilon,n+1}_{j,h}, v_h)
- b^*(U^{\epsilon,n}_j, u^{\epsilon,n}_{j,h}, v_h).
\end{gathered}
\end{equation}
Set $v_h = \phi^{\epsilon,n+1}_{j,h}$ and $q^h = p^{\epsilon,n+1}_{j, h}-\Tilde{q}$, then apply the polarization identity. We have
\begin{equation} 
\begin{gathered}
\frac{1}{2 \Delta t}( \|\phi^{\epsilon,n+1}_{j,h}\|^2 -\|\phi^{\epsilon,n}_{j,h}\|^2 + \|\phi^{\epsilon,n+1}_{j,h}- \phi^{\epsilon,n}_{j,h}\|^2)+\nu \|\nabla \phi^{\epsilon,n+1}_{j,h}\|^2 +\epsilon \|p^{\epsilon,n+1}_{j,h}-\Tilde{q}\|^2\\
=\frac{1}{\Delta t}(\eta^{\epsilon,n+1}_j-\eta^{\epsilon,n}_j,\phi^{\epsilon,n+1}_{j,h})+ \epsilon (p^{\epsilon,n+1}_{j,h}-\Tilde{q},p^{\epsilon,n+1}_{j}-\Tilde{q})+(r^{\epsilon,n+1}_j, \phi^{\epsilon,n+1}_{j,h})\\ 
+ b^*(u^{\epsilon, n+1}_j, u^{\epsilon,n+1}_j, \phi^{\epsilon,n+1}_{j,h}) - b^*(\langle u^{\epsilon}_h\rangle^n, u^{\epsilon,n+1}_{j,h}, \phi^{\epsilon,n+1}_{j,h}) - b^*(U^{\epsilon,n}_j, u^{\epsilon,n}_{j,h}, \phi^{\epsilon,n+1}_{j,h}).
\end{gathered}
\end{equation}
We bound the terms on the right--hand side.\\
$\frac{1}{\Delta t} (\eta^{\epsilon,n+1}_j-\eta^{\epsilon,n}_j,\phi^{\epsilon,n+1}_{j,h})$ term:
\begin{equation*}
\begin{gathered}
\frac{1}{\Delta t} (\eta^{\epsilon,n+1}_j-\eta^{\epsilon,n}_j,\phi^{\epsilon,n+1}_{j,h})
\leq \| \frac{\eta^{\epsilon,n+1}_j-\eta^{\epsilon,n}_j}{\Delta t}\|_{-1}\|\nabla \phi^{\epsilon,n+1}_{j,h}\|\\
\leq C(\nu)\|\frac{\eta^{\epsilon,n+1}_j-\eta^{\epsilon,n}_j}{\Delta t}\|^2_{-1} + \frac{\nu}{44} \|\nabla \phi^{\epsilon,n+1}_{j,h}\|^2\\
\leq C(\nu)\|\frac{\eta^{\epsilon,n+1}_j-\eta^{\epsilon,n}_j}{\Delta t}\|^2 + \frac{\nu}{44} \|\nabla \phi^{\epsilon,n+1}_{j,h}\|^2.
\end{gathered}
\end{equation*}
By the integral form of Taylor's theorem, we have
\begin{equation*}
\begin{gathered}
\eta^{\epsilon,n+1}_j=\eta^{\epsilon,n}_j + \int_{t_n}^{t_{n+1}} \eta^{\epsilon}_{j,t}\, ds.\\
\end{gathered}
\end{equation*}
Divided by $\Delta t $ on both sides, and take the $L^2$ norm on $D$,
\begin{equation*}
\begin{gathered}
\|\frac{\eta^{\epsilon,n+1}_j-\eta^{\epsilon,n}_j}{\Delta t}\|^2
=\int_{D} \left(\frac{1}{\Delta t} \int_{t_n}^{t_{n+1}} \eta^{\epsilon}_{j,t}\, ds\right)^2\, dx\\
\leq \frac{1}{(\Delta t)^2}\int_{D} \int_{t_n}^{t_{n+1}} 1\, ds \int_{t_n}^{t_{n+1}}|\eta^\epsilon_{j,t}|^2\, ds\, dx\\
\leq \frac{1}{\Delta t}\int_{D} \int_{t_n}^{t_{n+1}} |\eta^{\epsilon}_{j,t}|^2\, ds \, dx.
\end{gathered}
\end{equation*}
By Fubini's theorem, we have
\begin{equation*}
\begin{gathered}
\|\frac{\eta^{\epsilon,n+1}_j-\eta^{\epsilon,n}_j}{\Delta t}\|^2
\leq \frac{1}{\Delta t} \int_{t_n}^{t_{n+1}}\int_{D} |\eta^\epsilon_{j,t}|^2 \, dx\, ds
\leq \max_{t_n\leq t\leq t_{n+1}} \|\eta^{\epsilon}_{j,t}\|^2.
\end{gathered}
\end{equation*}
Thus, we have 
\begin{equation*}
\begin{gathered}
\frac{1}{\Delta t} (\eta^{\epsilon,n+1}_j-\eta^{\epsilon,n}_j,\phi^{\epsilon,n+1}_{j,h})
\leq C(\nu)\max_{t_n\leq t\leq t_{n+1}} \|\eta^{\epsilon}_{j,t}\|^2 + \frac{\nu}{44}\|\nabla \phi^{\epsilon,n+1}_{j,h}\|^2.
\end{gathered}
\end{equation*}
$(r^{\epsilon,n+1}_j, \phi^{\epsilon,n+1}_{j,h})$ term:
\begin{equation}
\begin{gathered}
(r^{\epsilon,n+1}_j, \phi^{\epsilon,n+1}_{j,h})
\leq \|r^{\epsilon,n+1}_j\|_{-1} \|\nabla  \phi^{\epsilon,n+1}_{j,h}\|\\
\leq C(\nu)\|r^{\epsilon,n+1}_j\|^2_{-1}+ \frac{\nu}{44}\|\nabla \phi^{\epsilon,n+1}_{j,h}\|^2 \\
\leq C(\nu)\|r^{\epsilon,n+1}_j\|^2+ \frac{\nu}{44}\|\nabla \phi^{\epsilon,n+1}_{j,h}\|^2.
\end{gathered}
\end{equation}
Recall $r^{\epsilon,n+1}_j = u^{\epsilon,n+1}_{j,t} -\frac{u^{\epsilon,n+1}_j - u^{\epsilon,n}_j}{\Delta t}$.
By the integral form of Taylor's theorem:
\begin{equation}
\begin{gathered}
u^{\epsilon,n}_j = u^{\epsilon,n+1}_{j} -\Delta t u^{\epsilon,n+1}_{j,t} -\int_{t_n}^{t_{n+1}} u^{\epsilon}_{j,tt} (t_n -s)\, ds, \\
r^{\epsilon,n+1}_j = \frac{1}{\Delta t} \int_{t_n}^{t_{n+1}} u^{\epsilon}_{j,tt}(s-t_n)\, ds.
\end{gathered}
\end{equation}
\begin{equation}
\begin{gathered}
\|r^{\epsilon,n+1}_j\|^2
= \int_{D} \left(\frac{1}{\Delta t} \int_{t_n}^{t_{n+1}} u^{\epsilon}_{j,tt}(s-t_n)\, ds\right)^2 \, dx\\
\leq \int_{D} \left(\int_{t_n}^{t_{n+1}} u^{\epsilon}_{j,tt}\, ds\right)^2 \, dx
\leq  \int_{D} \left(\int_{t_n}^{t_{n+1}} |u^{\epsilon}_{j,tt}|\, ds\right)^2 \, dx\\
\leq \int_{D} \int_{t_n}^{t_{n+1}} 1\, ds\int_{t_n}^{t_{n+1}} |u^{\epsilon}_{j,tt}|^2\, ds \, dx
=\Delta t \int_{D}\int_{t_n}^{t_{n+1}} |u^{\epsilon}_{j,tt}|^2\, \, ds \, dx.
\end{gathered}
\end{equation}
By Fubini's theorem, we have
\begin{equation*}
\|r^{\epsilon,n+1}_j\|^2 \leq \Delta t \int_{t_n}^{t_{n+1}} \int_{D} |u^{\epsilon}_{j,tt}|^2 \, dx \, ds.
\end{equation*}
Hence
\begin{equation}
\begin{gathered}
(r^{\epsilon,n+1}_j, \phi^{\epsilon, n+1}_{j,h})\leq C(\nu) (\Delta t)^2  \int_{t_n}^{t_{n+1}}\int_{D} |u^{\epsilon}_{j,tt}|^2 \, dx \, ds + \frac{\nu}{44} \|\nabla \phi^{\epsilon,n+1}_{j,h}\|^2.
\end{gathered}
\end{equation}
$\epsilon (p^{\epsilon,n+1}_{j,h}-\Tilde{q},p^{\epsilon,n+1}_{j}-\Tilde{q})$ term:
\begin{equation*}
\begin{gathered}
\epsilon (p^{\epsilon,n+1}_{j,h}-\Tilde{q},p^{\epsilon,n+1}_{j}-\Tilde{q})
\leq \frac{\epsilon}{2} \|p^{\epsilon,n+1}_{j,h}-\Tilde{q}\|^2 + \frac{\epsilon}{2} \|p^{\epsilon,n+1}_{j}-\Tilde{q}\|^2.
\end{gathered}
\end{equation*}
 \item Last we bound the trilinear forms, i.e. $b^*(\cdot, \cdot,\cdot)$. Denote 
 \begin{equation*}
 A:= b^*(u^{\epsilon,n+1}_j, u^{\epsilon,n+1}_j, \phi^{\epsilon,n+1}_{j,h})-b^*(\langle u^{\epsilon}_h\rangle^n, u^{\epsilon,n+1}_{j,h}, \phi^{\epsilon,n+1}_{j,h})
-b^*(U^{\epsilon,n}_j, u^{\epsilon,n}_{j,h}, \phi^{\epsilon,n+1}_{j,h}).
\end{equation*}
First, we add and subtract $b^*(u^{\epsilon,n}_{j,h}, u^{\epsilon,n+1}_{j,h}, \phi^{\epsilon, n+1}_{j,h})$, and add $b^*(u^{\epsilon,n}_{j,h}, \phi^{\epsilon, n+1}_{j,h}, \phi^{\epsilon, n+1}_{j,h})=0$. We have
\begin{equation*}
\begin{gathered}
A=b^*(u^{\epsilon,n+1}_j, u^{\epsilon,n+1}_j, \phi^{\epsilon, n+1}_{j,h})-b^*(u^{\epsilon,n}_{j,h}, u^{\epsilon,n+1}_{j,h}-\phi^{\epsilon, n+1}_{j,h}, \phi^{\epsilon, n+1}_{j,h})
+b^*(U^{\epsilon, n}_j, u^{\epsilon,n+1}_{j,h} -u^{\epsilon,n}_{j,h}, \phi^{\epsilon, n+1}_{j,h}).
\end{gathered}
\end{equation*}
Since $u^{\epsilon,n+1}_j - u^{\epsilon,n+1}_{j,h}=\eta^{\epsilon, n+1}_j-\phi^{\epsilon, n+1}_{j,h}$, $u^{\epsilon,n+1}_{j,h}-\phi^{\epsilon, n+1}_{j,h}=u^{\epsilon,n+1}_j-\eta^{\epsilon, n+1}_j$. We have
\begin{equation*}
\begin{gathered}
A=b^*(u^{\epsilon,n+1}_j, u^{\epsilon,n+1}_j, \phi^{\epsilon, n+1}_{j,h})-b^*(u^{\epsilon,n}_{j,h}, u^{\epsilon,n+1}_j-\eta^{\epsilon, n+1}_j, \phi^{\epsilon, n+1}_{j,h})+b^*(U^{\epsilon, n}_j, u^{\epsilon,n+1}_{j,h} -u^{\epsilon,n}_{j,h}, \phi^{\epsilon, n+1}_{j,h})\\
=b^*(u^{\epsilon,n+1}_j-u^{\epsilon,n}_{j,h}, u^{\epsilon,n+1}_j, \phi^{\epsilon, n+1}_{j,h}) + b^*(u^{\epsilon,n}_{j,h}, \eta^{\epsilon, n+1}_j, \phi^{\epsilon, n+1}_{j,h})+b^*(U^{\epsilon, n}_j, u^{\epsilon,n+1}_{j,h} -u^{\epsilon,n}_{j,h}, \phi^{\epsilon, n+1}_{j,h}).
\end{gathered}
\end{equation*}
We add and subtract $b^*(u^{\epsilon, n}_j, u^{n+1}_{j},\phi^{\epsilon, n+1}_{j,h})$,
\begin{equation}
\begin{gathered}
A=b^*(u^{\epsilon,n+1}_j-u^{\epsilon, n}_j, u^{\epsilon, n+1}_{j},\phi^{\epsilon, n+1}_{j,h})+b^*(u^{\epsilon, n}_j-u^{\epsilon,n}_{j,h}, u^{\epsilon,n+1}_j,\phi^{\epsilon, n+1}_{j,h})\\
+b^*(u^{\epsilon,n}_{j,h},\eta^{\epsilon, n+1}_j,\phi^{\epsilon, n+1}_{j,h})
+b^*(U^{\epsilon, n}_j, u^{\epsilon,n+1}_{j,h} -u^{\epsilon,n}_{j,h}, \phi^{\epsilon, n+1}_{j,h}).
\end{gathered}
\end{equation}
Denote 
\begin{equation*}
\begin{gathered}
A_1:=b^*(u^{\epsilon,n+1}_j-u^{\epsilon, n}_j, u^{\epsilon, n+1}_{j},\phi^{\epsilon, n+1}_{j,h}),\
A_2:= b^*(u^{\epsilon, n}_j-u^{\epsilon,n}_{j,h}, u^{\epsilon,n+1}_j,\phi^{\epsilon, n+1}_{j,h}),\\
A_3: = b^*(u^{\epsilon,n}_{j,h},\eta^{\epsilon, n+1}_j,\phi^{\epsilon, n+1}_{j,h}),\
A_4:=b^*(U^{\epsilon, n}_j, u^{\epsilon,n+1}_{j,h} -u^{\epsilon,n}_{j,h}, \phi^{\epsilon, n+1}_{j,h}).
\end{gathered}
\end{equation*}
We estimate $A_i$, where $i= 1,\ldots, 4$, as follows.
First, we bound $A_1$.
\begin{equation}
\begin{gathered}
A_1\leq C\|\nabla (u^{\epsilon,n+1}_j - u^{\epsilon, n}_j)\| \|\nabla u^{\epsilon,n+1}_j\|\|\nabla \phi^{\epsilon, n+1}_{j,h}\|\\
\leq \frac{\nu}{44}\|\nabla \phi^{\epsilon, n+1}_{j,h}\|^2 + C(\nu) \|\nabla (u^{\epsilon,n+1}_j -u^{\epsilon, n}_j)\|^2 \|\nabla u^{\epsilon,n+1}_j\|^2\\
\leq \frac{\nu}{44}\|\nabla \phi^{\epsilon, n+1}_{j,h}\|^2+ C(\nu)\Delta t \left( \int_{t_n}^{t_{n+1}} \|\nabla u^{\epsilon}_{j,t}\|^2\, dt \right)\|\nabla u^{\epsilon,n+1}_j\|^2\\
\leq \frac{\nu}{44}\|\nabla \phi^{\epsilon, n+1}_{j,h}\|^2+C(\nu) (\Delta t)^2 \max_{t_n\leq t\leq t_{n+1}} \|\nabla u^{\epsilon}_{j,t}\|^2 \|\nabla u^{\epsilon,n+1}_j\|^2.
\end{gathered}
\end{equation}
We bound $A_2$.
\begin{equation*}
A_2= b^*(\eta^{\epsilon, n}_j, u^{\epsilon,n+1}_j, \phi^{\epsilon, n+1}_{j,h})-b^*(\phi^{\epsilon, n}_{j,h}, u^{\epsilon,n+1}_j, \phi^{\epsilon, n+1}_{j,h}). 
\end{equation*}
\begin{equation}
\begin{gathered}
b^*(\eta^{\epsilon, n}_j, u^{\epsilon,n+1}_j,\phi^{\epsilon, n+1}_{j,h})
\leq C\|\eta^{\epsilon, n}_j\| \|\nabla^{n+1}_j\|\|\nabla \phi^{\epsilon, n+1}_{j,h}\|\\
\leq \frac{\nu}{44} \|\nabla \phi^{\epsilon, n+1}_{j,h}\|^2 + C(\nu) \|\nabla u^{n+1}_{j}\|^2 \|\nabla \eta^{\epsilon, n}_j\|^2.
\end{gathered}
\end{equation}
\begin{equation}
\begin{gathered}
-b^*(\phi^{\epsilon, n}_{j,h},u^{\epsilon,n+1}_j,\phi^{\epsilon, n+1}_{j,h})
\leq C \sqrt{\|\nabla \phi^{\epsilon, n}_{j,h}\|\|\phi^{\epsilon, n}_{j,h}\|}\|\nabla u^{\epsilon,n+1}_j\|\|\nabla \phi^{\epsilon, n+1}_{j,h}\|\\
\leq \frac{\nu}{44}\|\nabla \phi^{\epsilon, n+1}_{j,h}\|^2 + C(\nu)\|\nabla \phi^{\epsilon, n}_{j,h}\|\|\phi^{\epsilon, n}_{j,h}\|\|\nabla u^{\epsilon,n+1}_j\|^2\\
\leq \frac{\nu}{44}\|\nabla \phi^{\epsilon, n+1}_{j,h}\|^2 + \frac{\nu}{4}\|\nabla \phi^{\epsilon, n}_{j,h}\|^2 + C(\nu)\|\phi^{\epsilon, n}_{j,h}\|^2\|\nabla u^{\epsilon,n+1}_j\|^4.
\end{gathered}
\end{equation}
Now we bound $A_3$.
\begin{equation*}
A_3=-b^*(\eta^{\epsilon, n}_j,\eta^{n+1}_{j},\phi^{\epsilon, n+1}_{j,h})+b^*(\phi^{\epsilon, n}_{j,h},\eta^{n+1}_{j},\phi^{\epsilon, n+1}_{j,h})+b^*(u^{\epsilon, n}_j,\eta^{n+1}_{j},\phi^{\epsilon, n+1}_{j,h}).
\end{equation*}
\begin{equation}
\begin{gathered}
-b^*(\eta^{\epsilon, n}_j, \eta^{\epsilon, n+1}_j,\phi^{\epsilon, n+1}_{j,h})
\leq C\|\nabla \eta^{\epsilon, n}_j\|\|\nabla \eta^{\epsilon, n+1}_j\|\|\nabla \phi^{\epsilon, n+1}_{j,h}\|\\
\leq \frac{\nu}{44}\|\nabla \phi^{\epsilon, n+1}_{j,h}\|^2 + C(\nu)\|\nabla \eta^{\epsilon, n}_j\|^2 \|\nabla \eta^{\epsilon, n+1}_j\|^2.\\
\end{gathered}
\end{equation}
\begin{equation}
\begin{gathered}
b^*(\phi^{\epsilon, n}_{j,h}, \eta^{\epsilon, n+1}_j,\phi^{\epsilon, n+1}_{j,h})
\leq \sqrt{\|\nabla \phi^{\epsilon, n}_{j,h}\|\|\phi^{\epsilon, n}_{j,h}\|}\|\nabla \eta^{\epsilon, n+1}_j\|\|\nabla \phi^{\epsilon, n+1}_{j,h}\|\\
\leq \frac{\nu}{44}\|\nabla \phi^{\epsilon, n+1}_{j,h}\|^2 + C(\nu)\|\nabla \phi^{\epsilon, n}_{j,h}\|\|\phi^{\epsilon, n}_{j,h}\|\|\nabla \eta^{\epsilon, n+1}_j\|^2\\
\leq \frac{\nu}{44}\|\nabla \phi^{\epsilon, n+1}_{j,h}\|^2+\frac{\nu}{4}\|\nabla \phi^{\epsilon, n}_{j,h}\|^2+ C(\nu)\|\nabla \eta^{\epsilon, n+1}_j\|^4 \|\phi^{\epsilon, n}_{j,h}\|^2.
\end{gathered}
\end{equation}
\begin{equation}
\begin{gathered}
b^*(u^{n}_{j}, \eta^{\epsilon, n+1}_j,\phi^{\epsilon, n+1}_{j,h})
\leq C\|\nabla u^{n}_{j}\|\|\nabla \eta^{\epsilon, n+1}_j\|\|\nabla \phi^{\epsilon, n+1}_{j,h}\|\\
\leq \frac{\nu}{44}\|\nabla \phi^{\epsilon, n+1}_{j,h}\|^2 + C(\nu) \|\nabla u^{\epsilon, n}_j\|^2 \|\nabla \eta^{\epsilon, n+1}_j\|^{2}.
\end{gathered}
\end{equation}
Last, we bound $A_4$.
\begin{equation*}
\begin{gathered}
A_4=b^*(U^{\epsilon, n}_j, u^{\epsilon,n+1}_j -u^{\epsilon, n}_j, \phi^{\epsilon, n+1}_{j,h})-b^*(U^{\epsilon, n}_j, \eta^{\epsilon, n+1}_j, \phi^{\epsilon, n+1}_{j,h})\\
+b(U^{\epsilon, n}_j,\eta^{\epsilon, n}_j,\phi^{\epsilon, n+1}_{j,h})
-b^*(U^{\epsilon, n}_j,\phi^{\epsilon, n}_{j,h},\phi^{\epsilon, n+1}_{j,h}).
\end{gathered}
\end{equation*}
\begin{equation}
\begin{gathered}
b^*(U^{\epsilon, n}_j, u^{\epsilon,n+1}_j - u^{\epsilon, n}_j,\phi^{\epsilon, n+1}_{j,h})
\leq C\|\nabla U^{\epsilon, n}_j\| \|\nabla(u^{\epsilon,n+1}_j - u^{\epsilon, n}_j)\| \|\nabla \phi^{\epsilon, n+1}_{j,h}\|\\
\leq \frac{\nu}{44}\|\nabla \phi^{\epsilon, n+1}_{j,h}\|^2 + C(\nu)\|\nabla U^{\epsilon, n}_j\|^2 \|\nabla (u^{\epsilon,n+1}_j -u^{\epsilon, n}_j)\|^2\\
\leq \frac{\nu}{44}\|\nabla \phi^{\epsilon, n+1}_{j,h}\|^2+ C(\nu)\Delta t \|\nabla U^{\epsilon, n}_j\|^2 (\int_{t_n}^{t_{n+1}} \|\nabla u^{\epsilon}_{j,t}\|^2 \, dt )\\
\leq \frac{\nu}{44}\|\nabla \phi^{\epsilon, n+1}_{j,h}\|^2 + C(\nu)(\Delta t)^2\|\nabla U^{\epsilon, n}_j\|^2 \max_{t_n\leq t\leq t_{n+1}} \|\nabla u^{\epsilon}_{j,t}\|^2.  
\end{gathered}
\end{equation}
\begin{equation}
\begin{gathered}
-b^*(U^{\epsilon, n}_j, \eta^{\epsilon, n+1}_j, \phi^{\epsilon, n+1}_{j,h})
\leq C\|\nabla U^{\epsilon, n}_j\|\|\nabla \eta^{\epsilon, n+1}_j\|\|\nabla \phi^{\epsilon, n+1}_{j,h}\|\\
\leq \frac{\nu}{44} \|\nabla \phi^{\epsilon, n+1}_{j,h}\|^2 + C(\nu)\|\nabla U^{\epsilon, n}_j\|^2\|\nabla \eta^{\epsilon, n+1}_j\|^2.
\end{gathered}
\end{equation}
\begin{equation}
\begin{gathered}
b^*(U^{\epsilon, n}_j, \eta^{\epsilon, n}_j, \phi^{\epsilon, n+1}_{j,h})
\leq C\|\nabla U^{\epsilon, n}_j\|\|\nabla \eta^{\epsilon, n}_j\|\|\nabla \phi^{\epsilon, n+1}_{j,h}\|\\
\leq \frac{\nu}{44} \|\nabla \phi^{\epsilon, n+1}_{j,h}\|^2 + C(\nu)\|\nabla U^{\epsilon, n}_j\|^2\|\nabla \eta^{\epsilon, n}_j\|^2.
\end{gathered}
\end{equation}
\begin{equation}
\begin{gathered}
-b^*(U^{\epsilon, n}_j, \phi^{\epsilon, n}_{j,h},\phi^{\epsilon, n+1}_{j,h})
=b^*(U^{\epsilon, n}_j, \phi^{\epsilon, n+1}_{j,h}-\phi^{\epsilon, n}_{j,h},\phi^{\epsilon, n+1}_{j,h}).
\end{gathered}
\end{equation}
Since $b^*(u,v,w)+b^*(u,w,v)=0$, we have
\begin{equation}
\begin{gathered}
-b^*(U^{\epsilon, n}_j, \phi^{\epsilon, n}_{j,h},\phi^{\epsilon, n+1}_{j,h})=b^*(U^{\epsilon, n}_j, \phi^{\epsilon,n+1}_{j,h}, \phi^{\epsilon, n+1}_{j,h}-\phi^{\epsilon, n}_{j,h} )\\
\leq C \|\nabla U^{\epsilon, n}_j\| \|\nabla \phi^{\epsilon, n+1}_{j,h}\| \sqrt{\| \nabla (\phi^{\epsilon, n+1}_{j,h}-\phi^{\epsilon, n}_{j,h})\| \|\phi^{\epsilon, n+1}_{j,h}-\phi^{\epsilon, n}_{j,h}\|} \\
\leq C \|\nabla U^{\epsilon, n}_j\|\|\nabla \phi^{\epsilon, n+1}_{j,h}\| \frac{1}{\sqrt{h}}\|\phi^{\epsilon, n+1}_{j,h}-\phi^{\epsilon, n}_{j,h}\|\\
\leq \frac{C\Delta t}{ h} \|\nabla U^{\epsilon, n}_j\|^2\|\nabla \phi^{\epsilon, n+1}_{j,h}\|^2 + \frac{1}{4\Delta t}\|\phi^{\epsilon, n+1}_{j,h}-\phi^{\epsilon, n}_{j,h}\|^2.
\end{gathered}
\end{equation}
Combine terms,
\begin{equation}\label{enery_estimate}
\begin{gathered}
\frac{1}{2 \Delta t}( \|\phi^{\epsilon,n+1}_{j,h}\|^2 -\|\phi^{\epsilon,n}_{j,h}\|^2 + \frac{1}{2}\|\phi^{\epsilon,n+1}_{j,h}- \phi^{\epsilon,n}_{j,h}\|^2)+\left(\frac{\nu}{4}- \frac{C\Delta t}{ h} \|\nabla U^{\epsilon, n}_j\|^2\right)\|\nabla \phi^{\epsilon,n+1}_{j,h}\|^2 \\
+ \frac{\nu}{2} (\|\nabla \phi^{\epsilon,n+1}_{j,h}\|^2-\|\nabla \phi^{\epsilon,n}_{j,h}\|^2)
+\frac{\epsilon}{2} \|p^{\epsilon,n+1}_{j,h}-\Tilde{q}\|^2 \\
\leq C(\nu)\max_{t_n\leq t\leq t_{n+1}} \|\eta^{\epsilon}_{j,t}\|^2
+C(\nu) (\Delta t)^2 \int_{t_n}^{t_{n+1}} \int_{D} |u^{\epsilon}_{j,tt}|^2 \, dx \, ds  + \frac{\epsilon}{2} \|p^{\epsilon,n+1}_{j}-\Tilde{q}\|^2 \\
+C(\nu) (\Delta t)^2 \max_{t_n\leq t\leq t_{n+1}} \|\nabla u^{\epsilon}_{j,t}\|^2 \|\nabla u^{\epsilon,n+1}_j\|^2+ C(\nu) \|\nabla u^{n+1}_{j}\|^2 \|\nabla \eta^{\epsilon, n}_j\|^2\\
+C(\nu)\|\phi^{\epsilon, n}_{j,h}\|^2\|\nabla u^{\epsilon,n+1}_j\|^4+ C(\nu)\|\nabla \eta^{\epsilon, n}_j\|^2 \|\nabla \eta^{\epsilon, n+1}_j\|^2
+C(\nu)\|\nabla \eta^{\epsilon, n+1}_j\|^4 \|\phi^{\epsilon, n}_{j,h}\|^2\\
+C(\nu) \|\nabla u^{\epsilon, n}_j\|^2 \|\nabla \eta^{\epsilon, n+1}_j\|^{2}
+ C(\nu)(\Delta t)^2\|\nabla U^{\epsilon, n}_j\|^2 \max_{t_n\leq t\leq t_{n+1}} \|\nabla u^{\epsilon}_{j,t}\|^2\\
+ C(\nu)\|\nabla U^{\epsilon, n}_j\|^2\|\nabla \eta^{\epsilon, n+1}_j\|^2
+ C(\nu)\|\nabla U^{\epsilon, n}_j\|^2\|\nabla \eta^{\epsilon, n}_j\|^2.
\end{gathered}
\end{equation}
By the CFL condition, we have 
\begin{equation*}\frac{\nu}{4}-\frac{C \Delta t}{h} \|\nabla U^{\epsilon,n}_j\|^2 \geq C_0 \nu > 0,
\end{equation*}
for some constant $C_0>0$.

Recall equation (\ref{enery_estimate}), multiply by $2 \Delta t$ and organize terms:
\begin{equation}\label{convergence_n}
\begin{gathered}
\|\phi^{\epsilon, n+1}_{j,h}\|^2 -\|\phi^{\epsilon,n}_{j,h}\|^2 + \frac{1}{2}\|\phi^{\epsilon,n+1}_{j,h}- \phi^{\epsilon,n}_{j,h}\|^2+ C_0\Delta t \nu \|\nabla \phi^{\epsilon,n+1}_{j,h}\|^2\\
+\Delta t \bigg\{\nu (\|\nabla \phi^{\epsilon,n+1}_{j,h}\|^2 -\|\nabla \phi^{\epsilon,n}_{j,h}\|^2)+\epsilon \|p^{\epsilon,n+1}_{j,h}-\Tilde{q}\|^2\bigg\}\\
\leq  \Delta t  \bigg\{C(\nu)\left(\|\nabla u^{\epsilon,n+1}_{j}\|^4
 +\|\nabla \eta^{\epsilon,n+1}_j\|^4\right)\|\phi^{\epsilon,n}_{j,h}\|^2\\
+C(\nu) \max_{t_n\leq t\leq t_{n+1}}\| \eta^{\epsilon}_{j,t}\|^2
+C(\nu) (\Delta t)^2 \int_{t_n}^{t_{n+1}} \int_{D} |u^{\epsilon}_{j,tt}|^2 \, dx \, ds +\epsilon\|p^{\epsilon,n+1}_{j}-\Tilde{q}\|^2 \\
+ C(\nu)(\|\nabla \eta^{\epsilon,n}_j\|^2+\|\nabla u^{\epsilon,n}_j\|^2+\|\nabla U^{\epsilon,n}_j\|^2)\|\nabla \eta^{\epsilon,n+1}_{j}\|^2\\
+C(\nu)(\|\nabla u^{\epsilon,n+1}_j\|^2 +\|\nabla U^{\epsilon,n}_j\|^2)\|\nabla \eta^{\epsilon,n}_j\|^2\\
+ C(\nu)(\Delta t)^2(\|\nabla u^{\epsilon,n+1}_j\|^2+\|\nabla U^{\epsilon,n}_j\|^2)  \max_{t_n\leq t\leq t_{n+1}} \|\nabla u^{\epsilon}_{j,t}\|^2
\bigg\}.
\end{gathered}
\end{equation}
Take the sum of equation (\ref{convergence_n}) from $n=0$ to $n=N-1$, we have
\begin{equation*}
\begin{gathered}
\|\phi^{\epsilon, N}_{j,h}\|^2 + \frac{1}{2}\sum_{n=0}^{N-1} \|\phi^{\epsilon,n+1}_{j,h}- \phi^{\epsilon,n}_{j,h}\|^2 + \Delta t\nu \|\nabla \phi^{\epsilon,N}_{j,h}\|^2
+C_0\sum_{n=0}^{N-1}\Delta t \nu \|\nabla \phi^{\epsilon,n+1}_{j,h}\|^2\\
+ \Delta t \sum_{n=0}^{N-1} \epsilon \|p^{\epsilon,n+1}_{j,h}-\Tilde{q}\|^2
 \leq \|\phi^{\epsilon, 0}_{j,h}\|^2
 +\Delta t \nu \|\nabla \phi^{\epsilon,0}_{j,h}\|^2\\
+\Delta t \bigg\{\nu (\|\nabla \phi^{\epsilon,n+1}_{j,h}\|^2 -\|\nabla \phi^{\epsilon,n}_{j,h}\|^2)+\epsilon \|p^{\epsilon,n+1}_{j,h}-\Tilde{q}\|^2\bigg\}\\
\leq  \sum_{n=0}^{N-1}\Delta t  \bigg\{C(\nu)\left(\|\nabla u^{\epsilon,n+1}_{j}\|^4
 +\|\nabla\eta^{\epsilon,n+1}_j\|^4\right)\|\phi^{\epsilon,n}_{j,h}\|^2\\
+C(\nu) \max_{t_n\leq t\leq t_{n+1}}\| \eta^{\epsilon}_{j,t}\|^2
+C(\nu) (\Delta t)^2 \int_{t_n}^{t_{n+1}} \int_{D} |u^{\epsilon}_{j,tt}|^2 \, dx \, ds\\ +\epsilon\|p^{\epsilon,n+1}_{j}-\Tilde{q}\|^2
+ C(\nu)(\|\nabla \eta^{\epsilon,n}_j\|^2+\|\nabla u^{\epsilon,n}_j\|^2+\|\nabla U^{\epsilon,n}_j\|^2)\|\nabla \eta^{\epsilon,n+1}_{j}\|^2\\
+C(\nu)(\|\nabla u^{\epsilon,n+1}_j\|^2 + \|\nabla U^{\epsilon,n}_j\|^2)\|\nabla \eta^{\epsilon,n}_j\|^2\\
+ C(\nu)(\Delta t)^2(\|\nabla u^{\epsilon,n+1}_j\|^2+\|\nabla U^{\epsilon,n}_j\|^2)  \max_{t_n\leq t\leq t_{n+1}} \|\nabla u^{\epsilon}_{j,t}\|^2
\bigg\}.
\end{gathered}
\end{equation*}
By Lemma \ref{discrete_gronwall}, we have 
\begin{equation*}
\begin{gathered}
\|\phi^{\epsilon, N}_{j,h}\|^2 + \frac{1}{2}\sum_{n=0}^{N-1} \|\phi^{\epsilon,n+1}_{j,h}- \phi^{\epsilon,n}_{j,h}\|^2 + \Delta t \nu \|\nabla \phi^{\epsilon,N}_{j,h}\|^2++C_0\sum_{n=0}^{N-1}\Delta t \nu \|\nabla \phi^{\epsilon,n+1}_{j,h}\|^2
\\+ \Delta t \sum_{n=0}^{N-1} \epsilon \|p^{\epsilon,n+1}_{j,h}-\Tilde{q}\|^2
\leq \exp\Biggl\{C(\nu)\Delta t \sum_{n=0}^{N-1} \left(\|\nabla u^{\epsilon,n+1}_{j}\|^4
 +\|\nabla \eta^{\epsilon,n+1}_j\|^4\right)\Biggl\}\\
 \Biggl\{\|\phi^{\epsilon, 0}_{j,h}\|^2+ \Delta t \nu\|\nabla \phi^{\epsilon,0}_{j,h}\|^2 +\Delta t\sum_{n=0}^{N-1}\Biggl(C(\nu) \max_{t_n\leq t\leq t_{n+1}}\| \eta^{\epsilon}_{j,t}\|^2\\
+C(\nu) (\Delta t)^2 \int_{t_n}^{t_{n+1}} \int_{D} |u^{\epsilon}_{j,tt}|^2 \, dx \, ds +\epsilon\|p^{\epsilon,n+1}_{j}-\Tilde{q}\|^2\\
+ C(\nu)(\|\nabla \eta^{\epsilon,n}_j\|^2+\|\nabla u^{\epsilon,n}_j\|^2+\|\nabla U^{\epsilon,n}_j\|^2)\|\nabla \eta^{\epsilon,n+1}_{j}\|^2\\
+C(\nu)(\|\nabla u^{\epsilon,n+1}_j\|^2 + \|\nabla U^{\epsilon,n}_j\|^2)\|\nabla \eta^{\epsilon,n}_j\|^2\\
+ C(\nu)(\Delta t)^2(\|\nabla u^{\epsilon,n+1}_j\|^2+\|\nabla U^{\epsilon,n}_j\|^2)  \max_{t_n\leq t\leq t_{n+1}} \|\nabla u^{\epsilon}_{j,t}\|^2
\Biggl) \Biggl\}.
\end{gathered}
\end{equation*}
By Proposition \ref{average_equality}, we can conclude that 
\begin{equation*}
\Delta t\sum_{n=0}^{N} \|\nabla U^{\epsilon,n}_j\|^2 < C.
\end{equation*}
By Proposition \ref{stokes-estimate}, we have
\begin{equation*}
\begin{gathered}
\|\phi^{\epsilon, N}_{j,h}\|^2 + \frac{1}{2}\sum_{n=0}^{N-1} \|\phi^{\epsilon,n+1}_{j,h}- \phi^{\epsilon,n}_{j,h}\|^2 + \Delta t \nu \|\nabla \phi^{\epsilon,N}_{j,h}\|^2 + \Delta t \sum_{n=0}^{N-1} \epsilon \|p^{\epsilon,n+1}_{j,h}-\Tilde{q}\|^2\\
+C_0 \sum_{n=0}^{N-1}\Delta t \nu \|\nabla \phi^{\epsilon,n+1}_{j,h}\|^2 \leq \exp(\alpha )\Biggl\{\|\phi^{\epsilon, 0}_{j,h}\|^2+ \Delta t \nu \|\nabla \phi^{\epsilon,0}_{j,h}\|^2\\
+ C(\nu)T\left(\inf_{v_h\in X^h} \||\nabla (u^{\epsilon}_j-v_h)_t|\|^2_{\infty,0} +\inf_{q^h \in Q^h}\||(p^{\epsilon}_j -q^h)_t|\|^2_{\infty,0} \right)+ (\Delta t)^3 C(\nu) \||u^{\epsilon}_{j,tt}|\|^2_{\infty, 0}\\
+ \epsilon \Delta t C(\nu, \beta^h) \left(\inf_{v_h \in X^h}\||\nabla (u^{\epsilon}_j- v_h)_t|\|^2_{2,0}  + \inf_{q^h \in Q^h} \||(p^{\epsilon}_j - q^h)_t|\|^2_{2,0} \right)\\
 + C(\nu) \left(\inf_{v_h\in X^h}\||\nabla (u^{\epsilon}_j-v_h)|\|^2_{\infty,0} +  \inf_{q^h \in Q^h}\||p^{\epsilon}_j -q^h|\|^2_{\infty,0} +\||\nabla u^{\epsilon}_j|\|^2_{\infty,0}+CT\right)\\
 \left( \inf_{v_h\in X^h} \||\nabla (u^{\epsilon}_j-v_h)|\|^2_{2,0} + \inf_{q^h \in Q^h} \||p^{\epsilon}_j-q^h|\|^2_{2,0}\right)\\
 + C(\nu) (\Delta t)^2 \left(\Delta t \|\nabla u^{\epsilon}_j\|^2_{2, 0}+ C \right)  \||\nabla u^{\epsilon}_{j,t}|\|^2_{\infty,0}
\Biggl\},
\end{gathered}
\end{equation*}
where 
\begin{equation}
\alpha = C(\nu)\Delta t \sum_{n=0}^{N-1} \|\nabla u^{\epsilon,n+1}_{j}\|^4.
\end{equation}
Apply interpolation inequalities in equation (\ref{interpolation_inequality}),
\begin{equation*}
\begin{gathered}
\|\phi^{\epsilon, N}_{j,h}\|^2 + \frac{1}{2}\sum_{n=0}^{N-1} \|\phi^{\epsilon,n+1}_{j,h}- \phi^{\epsilon,n}_{j,h}\|^2 + \Delta t \nu \|\nabla \phi^{\epsilon,N}_{j,h}\|^2 + \Delta t \sum_{n=0}^{N-1} \epsilon \|p^{\epsilon,n+1}_{j,h}-\Tilde{q}\|^2\\
+C_0\Delta t \sum_{n=0}^{N-1}\nu\|\nabla \phi^{\epsilon,n+1}_{j,h}\|^2 \leq \exp (\alpha )\Biggl\{\|\phi^{\epsilon, 0}_{j,h}\|^2+ \Delta t \nu \|\nabla \phi^{\epsilon,0}_{j,h}\|^2\\
+ h^{2m}C(\nu) T\left( \||u^{\epsilon}_{j,t}|\|^2_{\infty,0} + \frac{1}{\nu^2} \||p^{\epsilon}_{j,t}|\|^2_{\infty,0} \right)+ (\Delta t)^3C(\nu) \||u^{\epsilon}_{j,tt}|\|^2_{\infty, 0}\\
+ h^{2m}\epsilon \Delta t C(\nu, \beta^h)\left(\||u^{\epsilon}_{j,t}|\|^2_{2,0}  + \||p^{\epsilon}_{j,t}|\|^2_{2,0} \right)\\
+ h^{2m}C(\nu) T\left(  \|| (u^{\epsilon}_j|\|^2_{2,0} +\frac{1}{\nu^2} \||p^{\epsilon}_j|\|^2_{2,0}\right)+ C(\nu) (\Delta t)^2 \||\nabla u^{\epsilon}_{j,t}|\|^2_{\infty,0}
\Biggl\}
\end{gathered}
\end{equation*}
Recall that $e^{\epsilon,n}_{j}= \eta^{\epsilon,n}_j -\phi^{\epsilon,n}_{j,h}$. Using the triangle inequality, we have
\begin{equation*}
\begin{gathered}
\|e^{\epsilon, N}_{j}\|^2 + \frac{1}{2}\sum_{n=0}^{N-1} \|e^{\epsilon,n+1}_{j}- e^{\epsilon,n}_{j}\|^2 + \Delta t \nu \|\nabla e^{\epsilon,N}_{j}\|^2
+C_0\Delta t \sum_{n=0}^{N-1}\nu\|\nabla e^{\epsilon,n+1}_{j}\|^2\\
\leq \|\phi^{\epsilon, N}_{j,h}\|^2 + \frac{1}{2}\sum_{n=0}^{N-1} \|\phi^{\epsilon,n+1}_{j,h}- \phi^{\epsilon,n}_{j,h}\|^2 + \Delta t \nu \|\nabla \phi^{\epsilon,N}_{j,h}\|^2
+C_0\Delta t \sum_{n=0}^{N-1}\nu\|\nabla \phi^{\epsilon,n+1}_{j,h}\|^2\\
+\|\eta^{\epsilon, N}_{j}\|^2 + \frac{1}{2}\sum_{n=0}^{N-1} \|\eta^{\epsilon,n+1}_{j}- \eta^{\epsilon,n}_{j}\|^2 + \Delta t \nu \|\nabla \eta^{\epsilon,N}_{j}\|^2
+C_0\Delta t \sum_{n=0}^{N-1}\nu\|\nabla \eta^{\epsilon,n+1}_{j}\|^2.
\end{gathered}
\end{equation*}
We complete the proof using the previous bounds for the $\eta^\epsilon_{j}$ terms.
\end{proof}
Combining Theorem \ref{convergence_ensemble_penalty} with the result of Shen \cite{shen1995error}, Theorem $4.1$, p. 395, and applying the triangle inequality,
\begin{equation*}
\|u_j(t_n)-u^{\epsilon,n}_{j,h}\|\leq \|u_j(t_n)-u^\epsilon_j(t_n)\|+\|u^\epsilon_j(t_n)-u^{\epsilon,n}_{j,h}\|.
\end{equation*}
We have the following corollaries.
\begin{corollary}\label{final-error}
Assume the regular solutions, under the CFL condition in equation (\ref{stability_condition}, we have the following optimal estimates:
\begin{equation*}
\max_{t_n} \| u_j(t_n) - u^{\epsilon,n}_{j,h}\|^2 +\Delta t \sum_{n=1}^{N} \|\nabla ( u_j(t_n) - u^{\epsilon,n}_{j,h})\|^2 \leq C(u_j,\nu, T) (\epsilon +\Delta t + h^m)^2.
\end{equation*}
\end{corollary}
\begin{corollary} The error between the average of true solution and the average of penalized finite element approximations is
\begin{equation*}
\| \langle u_{t_n}\rangle- \langle u^{\epsilon}_h\rangle^n\|^2
\leq C(u_1,\ldots, u_J, \nu, T) (\epsilon +\Delta t + h^m)^2.
\end{equation*}
\end{corollary}
\begin{proof}
\begin{equation*}
\begin{gathered}
   \| \langle u(t_n)\rangle- \langle u^{\epsilon}_h\rangle^n\|^2=\|\frac{1}{J} \sum_{j=1}^J (u_j-u^{\epsilon,n}_{j,h})\|^2
   =\left(\frac{1}{J}\right)^2 \|\sum_{j=1}^J (u_j-u^{\epsilon,n}_{j,h})\|^2.
   \end{gathered}
   \end{equation*}
By the Cauchy Schwarz inequality,
\begin{equation*}
\|\sum_{j=1}^J (u_j-u^{\epsilon,n}_{j,h})\|^2\leq J \sum_{j=1}^J \|u_j-u^{\epsilon,n}_{j,h}\|^2.
\end{equation*}
By Corollary \ref{final-error},
\begin{equation*}
\sum_{j=1}^J \|u_j-u^{\epsilon,n}_{j,h}\|^2\leq J C(u_1,\ldots, u_J, \nu, T) (\epsilon +\Delta +h^m).
\end{equation*}
Thus,
\begin{equation*}
\sum_{j=1}^J \|u_j-u^{\epsilon,n}_{j,h}\|^2\leq J^2 C(u_1,\ldots, u_J, \nu, T) (\epsilon +\Delta +h^m).
\end{equation*}
Hence, we have
\begin{equation*}
   \begin{gathered}
      \| \langle u(t_n)\rangle- \langle u^{\epsilon}_h\rangle^n\|^2 \leq C(u_1,\ldots, u_J, \nu, T) (\epsilon +\Delta +h^m).
   \end{gathered}
\end{equation*}
\end{proof}
\section{Ensemble--based Monte Carlo forecasting}\label{sec: monte-carlo}
We consider the NSE with random body forces and initial conditions. We find random functions $u: \Omega \times \bar{D} \times [0,T] \to \mathbb{R}^d$, and $p: \Omega \times \bar{D} \times [0,T] \to \mathbb{R}$ satisfy
\begin{equation}
\begin{gathered}
\frac{\partial u}{\partial t} +u\cdot \nabla u - \nu \Delta u + \nabla p = f(\omega, x,t),\\
\nabla \cdot u =0.
\end{gathered}
\end{equation}
 We choose a set of random samples for the random body force $f_j \equiv f(\omega_j, \cdot, \cdot)$, initial condition $u^0_j \equiv u^0(\omega_j,\cdot, \cdot)$ for $j=1,\ldots, J$. Note that the corresponding solutions $u(\omega_j,\cdot, \cdot)$ are independent, identically distributed (i.i.d). 

The penalty--based ensemble Monte Carlo is defined as follows. Denote $u^{\epsilon,n}_{j,h} = u^\epsilon_h(\omega_j, x, t_n)$ and $p^{\epsilon,n}_{j,h}=p^\epsilon(\omega_j, x, t_n)$. For the $j_{th}$ ensemble member and for $0 \leq n\leq N-1$, find $(u^{\epsilon,n+1}_{j,h},p^{\epsilon,n+1}_{j,h})  \in (X^h, Q^h)$ satisfying: 
\begin{equation}\label{MC-weak_form_method}
\begin{gathered}
\frac{1}{\Delta t_n}(u^{\epsilon, n+1}_{j,h}-u^{\epsilon,n}_{j,h}, v_h) + b^*(\langle u^\epsilon_h\rangle^n, u^{\epsilon,n+1}_{j,h},v_h)+ b^*(u^{\epsilon,n}_{j,h} -\langle u^\epsilon_h\rangle^n, u^{\epsilon,n}_{j,h},v_h)\\
+\nu (\nabla u^{\epsilon,n+1}_{j,h}, \nabla v_h) -(p^{\epsilon,n+1}_{j,h}, \nabla \cdot v_h) +(q^h, \nabla \cdot u^{\epsilon, n+1}_{j,h}) + \epsilon (p^{\epsilon, n+1}_{j,h}, q^h)= (f_j^{n+1},v_h),
\end{gathered}
\end{equation}
for all $(v_h, q^h)\in (X^h,Q^h)$.

We approximate $E[u]$ by the sample average of the penalized NSE $\frac{1}{J} \sum_{j=1}^J u^\epsilon_h(\omega_j, \cdot,\cdot)$.

Theorem \ref{Stability of BEFE-Ensemble penalty} together with the property of expectation leads to the following stability analysis for the finite element solution $u^{\epsilon,n}_{j,h}$.
\begin{theorem}
\label{Stability_MC_Ensemble_penalty} Suppose the following timestep condition holds:
\begin{equation}
C\frac{\Delta t}{\nu h} E[\|\nabla U^{\epsilon,n}_j\|^2]\leq 1, j=1,\ldots, J.
\end{equation}
Then for any $N\geq 1$:
\begin{equation}
\begin{gathered}
\frac{1}{2} E[\|u^{\epsilon,N}_{j,h}\|^2] + \frac{1}{4}\sum_{n=0}^{ N-1}E[\|u^{\epsilon,n+1}_{j,h}- u^{\epsilon,n}_{j,h}\|^2] + \frac{\nu \Delta t}{4} E[\|\nabla u^{\epsilon,N}_{j,h}\|^2]\\
+ \frac{\Delta t}{\epsilon}\sum_{n=0}^{N-1} E[\|P_{Q^h}\nabla \cdot u^{\epsilon,n+1}_{j,n}\|^2] 
+\frac{\nu \Delta t}{4}\sum_{n=0}^{N-1} E[\|\nabla u^{\epsilon,n+1}_{j,h}\|^2 ]\\
\leq \frac{\Delta t}{2 \nu}\sum_{n=0}^{N-1} E[\|f^{n+1}_{j,h}\|^2_{-1}]+ \frac{1}{2} E[\|u^0_{j,h}\|^2 ]+\frac{\nu \Delta t}{4} E[\|\nabla u^0_{j,h}\|^2 ]
\end{gathered}
\end{equation}
\end{theorem}
The fully discrete penalty--based ensemble Monte Carlo approximation is defined to be 
\begin{equation*}
\Psi^n_{h}=\frac{1}{J}\sum_{j=1}^J u^{\epsilon,n}_{j,h}.
\end{equation*} 
We estimate $E[u^\epsilon(t_n)] -\Psi^n_h$ in averaged norms. We write
\begin{equation*}
\begin{gathered}
E[u^\epsilon(t_n)]- \Psi^n_h = (E[u^\epsilon(t_n)]-E[u^{\epsilon,n}_{j,h}]) + (E[u^{\epsilon,n}_{j,h}]-\Psi^n_h).
\end{gathered}
\end{equation*}
Since $u^\epsilon_j$ are i.i.d, $E[u^\epsilon(t_n)]=E[u^\epsilon_j(t_n)]$. Thus,
\begin{equation*}
\begin{gathered}
E[u^\epsilon(t_n)]- \Psi^n_h 
=\Gamma^n_h + \Gamma^n_S,
\end{gathered}
\end{equation*}
where $\Gamma^n_h = E[u_j(t_n)]-E[u^{\epsilon,n}_{j,h}]$ is the discretization error, and $\Gamma^n_S=E[u^{\epsilon,n}_{j,h}]-\Psi^n_h$ is the statistical error controls by the ensemble size. 
\begin{theorem} Assume the condition in equation (\ref{stability_condition}) holds for all $n$,
\begin{equation}
C\frac{\Delta t}{\nu h} E[\|\nabla U^{\epsilon,n}_j\|^2]\leq 1, j=1,\ldots, J,
\end{equation}
then there are positive constant $C$ and $C_0$ independent of the $h$ and $\Delta t$ such that
\begin{equation}
\begin{gathered}
E[\|e^{\epsilon, N}_{j,h}\|^2] + \frac{1}{2}\sum_{n=0}^{N-1} E[\|e^{\epsilon,n+1}_{j,h}- e^{\epsilon,n}_{j,h}\|^2] + \Delta t \nu E[\|\nabla e^{\epsilon,N}_{j,h}\|^2]\\
+ C_0\Delta t \sum_{n=0}^{N-1}\nu E[\|\nabla e^{\epsilon,n+1}_{j,h}\|^2] \leq \exp (\alpha )\Biggl\{E[\|e^{\epsilon, 0}_{j,h}\|^2]+ \Delta t \nu E[\|\nabla e^{\epsilon,0}_{j,h}\|^2]\\
+ h^{2m}C(\nu) T\left( E[\||u^{\epsilon}_{j,t}|\|^2_{\infty,0}] + \frac{1}{\nu^2} E[\||p^{\epsilon}_{j,t}|\|^2_{\infty,0}] \right)+ (\Delta t)^3C(\nu) E[\||u^{\epsilon}_{j,tt}|\|^2_{\infty, 0}]\\
+ h^{2m}\epsilon \Delta t C(\nu, \beta^h)\left(E[\||u^{\epsilon}_{j,t}|\|^2_{2,0}]  + E[\||p^{\epsilon}_{j,t}|\|^2_{2,0}] \right)\\
+ h^{2m}C(\nu) T\left( E[\|| u^{\epsilon}_j|\|^2_{2,0}] +\frac{1}{\nu^2} E[\||p^{\epsilon}_j|\|^2_{2,0}]\right)+ C(\nu) (\Delta t)^2 E[\||\nabla u^{\epsilon}_{j,t}|\|^2_{\infty,0}]
\Biggl\},
\end{gathered}
\end{equation}
where
\begin{equation*}
\alpha = C(\nu) \Delta t \sum_{n=0}^{N-1} E[\|\nabla u^{\epsilon, n+1}\|^4]. 
\end{equation*}
\end{theorem}
\begin{proof}
The conclusion follows Theorem \ref{convergence_ensemble_penalty} after applying the expectation on equation (\ref{convergence_eqn}).
\end{proof}
\begin{theorem} 
Consider the method in equation (\ref{penalty_ensemble}), assume that $\forall n$,
\begin{equation}
C\frac{\Delta t}{\nu h} E[\|\nabla U^{\epsilon,n}_j\|^2]\leq 1, \  j=1,\ldots, J.
\end{equation}
Then for any $N\geq 1$:
\begin{equation}\label{statistical-bound}
\begin{gathered}
\frac{1}{2} E[\|\Gamma^N_S\|^2] + \frac{1}{4}\sum_{n=0}^{ N-1}E[\|\Gamma^{n+1}_S- \Gamma^{n}_S\|^2] + \frac{\nu \Delta t}{4} E[\|\nabla \Gamma^{N}_S\|^2]\\
+ \frac{\Delta t}{\epsilon}\sum_{n=0}^{N-1} E[\|P_{Q^h}\nabla \cdot \Gamma^{n+1}_S\|^2] 
+\frac{\nu \Delta t}{4}\sum_{n=0}^{N-1} E[\|\nabla \Gamma^{n+1}_S\|^2 ]\\
\leq \frac{1}{J}\bigg\{\frac{\Delta t}{2 \nu}\sum_{n=0}^{N-1} E[\|f^{n+1}_{j,h}\|^2_{-1}]+ \frac{1}{2} E[\|u^0_{j,h}\|^2 ]+\frac{\nu \Delta t}{4} E[\|\nabla u^0_{j,h}\|^2 ]\bigg\}.
\end{gathered}
\end{equation}
\end{theorem}
\begin{proof}
Herein, we present the estimate $E[\|\nabla \Gamma^n_S\|^2]$. Define $\langle u^{\epsilon,n}_{j,h},u^{\epsilon,n}_{j,h}\rangle:=(\nabla u^{\epsilon,n}_{j,h}, \nabla u^{\epsilon,n}_{j,h})$.
\begin{equation*}
\begin{gathered}
E[\|\nabla \Gamma^n_S\|^2]= E\left[\langle \frac{1}{J} \sum_{i=1}^J (E[u^{\epsilon,n}_{i,h}]-u^{\epsilon, n}_{i,h}), \frac{1}{J} \sum_{i=1}^J (E[u^{\epsilon,n}_{j,h}]-u^{\epsilon, n}_{i,h})\rangle\right] \\
=\frac{1}{J^2} \sum_{i=1}^J \sum_{j=1}^J E[\langle E[u^{\epsilon,n}_{j,h}]-u^{\epsilon, n}_{j,h}, E[u^{\epsilon,n}_{j,h}]-u^{\epsilon, n}_{j,h}\rangle]\\
=\frac{1}{J^2} \sum_{j=1}^J E[\langle E[u^{\epsilon,n}_{j,h}]-u^{\epsilon, n}_{j,h}, E[u^{\epsilon,n}_{j,h}]-u^{\epsilon, n}_{j,h}\rangle].
\end{gathered}
\end{equation*}
The last equality is due to the fact $u^{\epsilon,n}_{j,h}$ for $j=1,\ldots, J$ are i.i.d., and when $i\neq j$, the expectation  of $\langle E[u^{\epsilon,n}_{j,h}]-u^{\epsilon, n}_{j,h}, E[u^{\epsilon,n}_{i,h}]-u^{\epsilon, n}_{i,h}\rangle$ is zero. We now expand the quantity $\langle E[u^{\epsilon,n}_{j,h}]-u^{\epsilon, n}_{j,h}, E[u^{\epsilon,n}_{j,h}]-u^{\epsilon, n}_{j,h}\rangle$. Use the fact $E[u^{\epsilon,n}_{h}] = E[u^{\epsilon,n}_{j,h}]$ and $E[(u^{\epsilon,n}_{h})^2]=E[(u^{\epsilon,n}_{j,h})^2]$ to obtain
\begin{equation*}
\begin{gathered}
E[\|\nabla \Gamma^n_S\|^2] = -\frac{1}{J} \|\nabla E[u^{\epsilon,n}_{j,h}]\|^2 + \frac{1}{J} E[\|\nabla u^{\epsilon,n}_{j,h}\|^2]\\
\leq \frac{1}{J} E[\|\nabla u^{\epsilon,n}_{j,h}\|^2].
\end{gathered}
\end{equation*}
The other terms involving the $E[\|\Gamma^N_S\|^2], E[\|\nabla \Gamma^N_S\|]$ and $E[\|\Gamma^{n+1}_{S}-\Gamma^{n}_{S}\|]$ can be treated similarly.
\end{proof}
The statistical error from sampling is $\mathcal{O}(\frac{1}{J})$. Combining Theorem \ref{convergence_ensemble_penalty} with the result of Shen \cite{shen1995error}, Theorem $4.1$, p. 395, and using the triangle inequality, we will have the following corollary. 
\begin{corollary}
\begin{equation*}
\begin{gathered}
\max_{t_n} E[\| u_j(t_n) - u^{\epsilon,n}_{j,h}\|^2] +\Delta t \sum_{n=1}^{N} E[\|\nabla ( u_j(t_n) - u^{\epsilon,n}_{j,h})\|^2] \leq C(u_j,\nu, T) (\epsilon +\Delta t + h^m)^2 \\
+\frac{1}{J}\bigg\{\frac{\Delta t}{2 \nu}\sum_{n=0}^{N-1} E[\|f^{n+1}_{j,h}\|^2_{-1}]+ \frac{1}{2} E[\|u^0_{j,h}\|^2 ]+\frac{\nu \Delta t}{4} E[\|\nabla u^0_{j,h}\|^2 ]\bigg\}.
\end{gathered}
\end{equation*}    
\end{corollary}
\section{Numerical Experiments} \label{sec:numerical-tests} We present the results of three numerical tests to illustrate our theory. In the first test, we calculate the rates of convergence using exact solutions with an ensemble size of two. Then, we construct a chaotic Lagrangian flow on a cylinder with perturbed body forces. In the third test, we extend this algorithm with the Coriolis force for a larger ensemble size, considering the benchmark test problem of flow past a cylinder. In these tests, we calculate various flow statistics to evaluate the flow dynamics:
\begin{equation*}
\begin{gathered}
|\text{angular momentum}| := |\int_{D} \bar{x} \times \bar{u}\, d \bar{x}|,\\
\text{enstrophy} := \frac{1}{2} \nu \|\nabla \times \bar{u}\|^2,\\
\text{kinetic energy} := \frac{1}{2} \| \bar{u} \|^2,\\
\text{viscous dissipation rate} : = \nu \|\nabla u\|^2,\\
\text{numerical dissipation rate from backward Euler (BE)} := \frac{1}{\Delta t} (u_n-u_{n-1})^2, \\
\text{numerical dissipation rate from penalizing incompressibility} := \frac{1}{\epsilon} \|\nabla \cdot u \|^2.
\end{gathered}
\end{equation*}
We use a second--order polynomial to approximate the velocity field in the following tests. The unstructured mesh is generated by GMSH \cite{geuzaine2009gmsh}. 

\subsection{Test for accuracy from \cite{ccibik2024ramshaw}}
\label{subsec:2}
We verify the convergence rates for the method in equation (\ref{weak_form_method}) with the following test. In $D = (0,1)^2$, the exact solution is given by 
\begin{align*}
\begin{gathered}
    u(x,y,t) =  (\exp(t) \cos(y), \exp(t) \sin(x))^\top,\\
    p(x,y,t) = (x-y) (1+t).
\end{gathered}
\end{align*}
The body force $f$ is calculated by substituting $u$ and $p$ in the NSE. We impose the Dirichlet boundary conditions where $
u_h = u_{true}$ on the boundary. We perturb the initial conditions  as follows:
\begin{align*}
\begin{gathered}
u_j(x,y,0) = (1 + \delta_j) u(x,y,0), \text{ for } j=1 \text{ and } 2,
\end{gathered}
\end{align*} 
where $\delta_1 = 10^{-3}$ and $\delta_2 = -10^{-3}$.

We set the kinematic viscosity $\nu = 1$, the characteristic velocity of the flow $U=1$, the characteristic length $L=1$, and the Reynolds number $Re = \frac{UL}{\nu}$. To discretize the domain, we choose a sequence of mesh sizes $h=\frac{1}{g}$, see Tables \ref{tab:1} and \ref{tab:2}.  We set  $\Delta t = \frac{h}{10}$,  $\epsilon= \Delta t$, and $T=1$. We denote the error as $e(h) = C h^\beta$. We solve the convergence rate $\beta$ via
\begin{equation*}
\beta = \frac{\ln(e(h_1)/e(h_2))}{\ln(h_1/h_2)},
\end{equation*}
at two successive values of $h$. Tables \ref{tab:1} and \ref{tab:2} show that the rates of convergence of $u_1$ and $u_2$ are optimal, second order.
\begin{table}[H]
\caption{The rates of convergence for $u_1$.}
\label{tab:1}
\centering
\begin{tabular}{c | c | c | c | c }
g & $\max_{t_n} \| u_1(t_n) - u^{\epsilon,n}_{1,h}\|$ & rate & $\sqrt{\Delta t \sum_{n=1}^{N} \|\nabla ( u_1(t_n) - u^{\epsilon,n}_{1,h})\|^2} $ & rate \\
\hline
$(\frac{3}{2})^0 \cdot27$ & 0.00358 & -- & 0.01353 & -- \\
$(\frac{3}{2})^1\cdot 27$ & 0.00169 & 1.91 & 0.00639 & 1.91\\
$(\frac{3}{2})^2 \cdot 27$ & 0.00076 & 1.95 & 0.0029 & 1.95\\
$(\frac{3}{2})^3 \cdot 27$ & 0.00033 & 1.98 & 0.00127 & 1.98\\
$(\frac{3}{2})^4 \cdot 27$ & 0.00015 & 1.99 & 0.00057 & 1.99
\end{tabular}
\end{table}
\begin{table}[H]
\caption{The rates of convergence for $u_2$.}
\label{tab:2}
\centering
\begin{tabular}{c | c | c | c | c }
g& $\max_{t_n} \| u_2(t_n) - u^{\epsilon,n}_{2,h}\|$ & rate & $\sqrt{\Delta t \sum_{n=1}^{N} \|\nabla ( u_2(t_n) - u^{\epsilon,n}_{2,h})\|^2}$ & rate \\
\hline
$(\frac{3}{2})^0 \cdot27$ & 0.00356& -- & 0.01348&-- \\
$(\frac{3}{2})^1\cdot 27$ & 0.00168 &1.91 & 0.00636 & 1.91\\
$(\frac{3}{2})^2 \cdot 27$ & 0.00076& 1.95& 0.00288 & 1.95\\
$(\frac{3}{2})^3 \cdot 27$ &0.00033&1.98& 0.00126 & 1.98 \\
$(\frac{3}{2})^4 \cdot 27$ & 0.00015 & 1.99& 0.00057 & 1.98
\end{tabular}
\end{table}

\subsection{Two rotating small cylinders} \label{sec: victor-aref}
\begin{figure}
\centering
   \begin{subfigure}{0.45\linewidth}
   \centering
\includegraphics[width=0.7\linewidth]{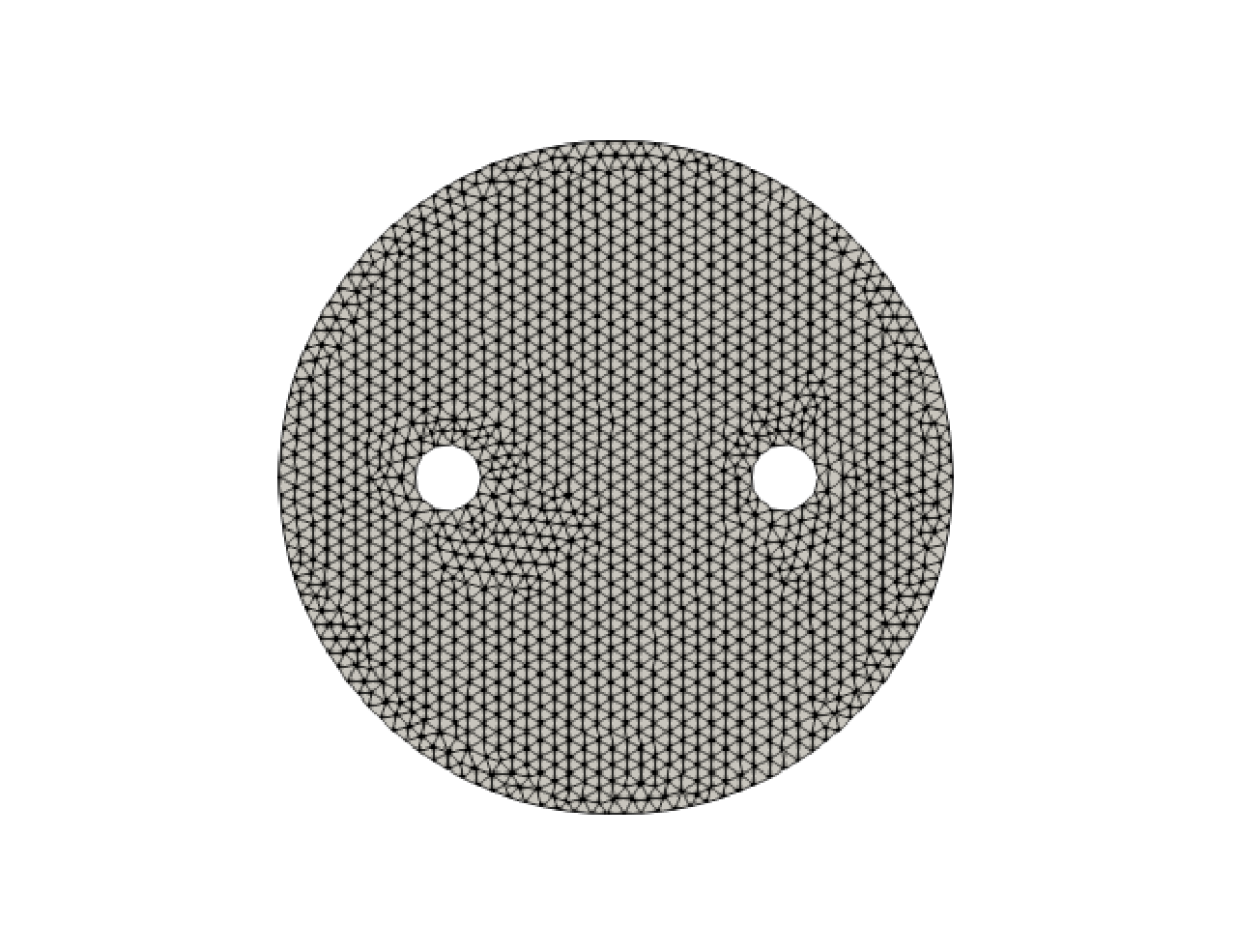}
   \caption{Mesh.}
   \label{fig: victor-mesh} 
\end{subfigure}
\hfill
\begin{subfigure}{0.45\linewidth}
   \centering
   \includegraphics[width=0.7\linewidth]{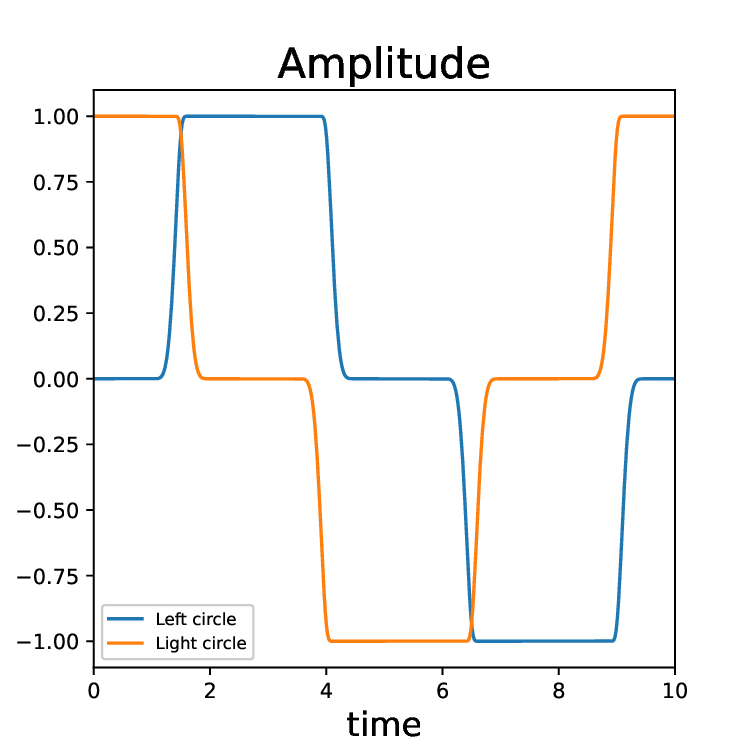}
   \caption{Amplitude.}
   \label{fig: amplitude}
\end{subfigure}
\caption{}
\label{fig: rotating-cylinders}
\end{figure}

We construct a simple 2D time--periodic flow that exhibits Lagrangian chaos, where the motion of fluid particles becomes chaotic, Aref \cite{aref1990chaotic}. Aref's blinking vortex flow is a model system to study chaotic advection and mixing in fluid flows, introduced by Aref \cite{ aref2020stirring, aref1983integrable}, and Aref and Balachandar \cite{aref1986chaotic}. The stirring was non--smooth over time, achieved using a point vortex. Herein, we use a cylinder with Dirichlet boundary conditions. The domain is a disk with two smaller obstacles inside (see Figure \ref{fig: victor-mesh}). We set the outer circle radius $r_0=1$, the left inner circle radius $r_1=0.1$, and the right inner circle radius $r_2=0.1$,  and $c=(c_1,c_2)=(\frac{1}{2}, 0)$. We define the domain:
\begin{equation*}
D = \{ (x,y): x^2 + y^2\leq r_0^2, (x+c_1)^2 + (y-c_2)^2 \geq r_1^2, \text{ and }  (x-c_1)^2 + (y-c_2)^2 \geq r_2^2 \}.
\end{equation*}
Dirichlet boundary conditions on the left and right circles rotate the flow. Figure \ref{fig: amplitude} shows the amplitude of the left and right circles. We have
\begin{equation*}
u(x,y) = 5 \text{ amplitude } (y, -x)^T \text{ on } \partial D.
\end{equation*}

Note that the outer circle remains
stationary. We chose mesh size $h= 0.05$, the final time $T = 10$, timestep $\Delta t = 0.001, \nu = 1/50$ and $\mathcal{R}e = 1/\nu$. The penalty parameter $\epsilon = \Delta t$. Flow is at rest at the beginning with exact boundary conditions. We perturbed the Dirichlet boundary conditions by the following:
\begin{equation*}
u_{1,2}(x,y) = (1 + \sigma_{1,2}) u(x,y) \text{ on } \partial D,
\end{equation*}
where $\sigma_1 = 0.01, \sigma_2 =-0.02$. We average the Dirichlet boundary for the ensemble members and write as
\begin{equation*}
u_{0}(x,y) = \frac{u_1(x,y) + u_2(x,y)}{2} \text{ on } \partial D.
\end{equation*}
We define the ensemble spread as follows:
\begin{equation*}
\text{ensemble spread}:=\frac{\|u_1-u_2\|}{\|u_{ave}\|}.
\end{equation*}

\begin{figure}
\centering
   \begin{subfigure}{0.45\linewidth}
   \centering
\includegraphics[width=0.7\linewidth]{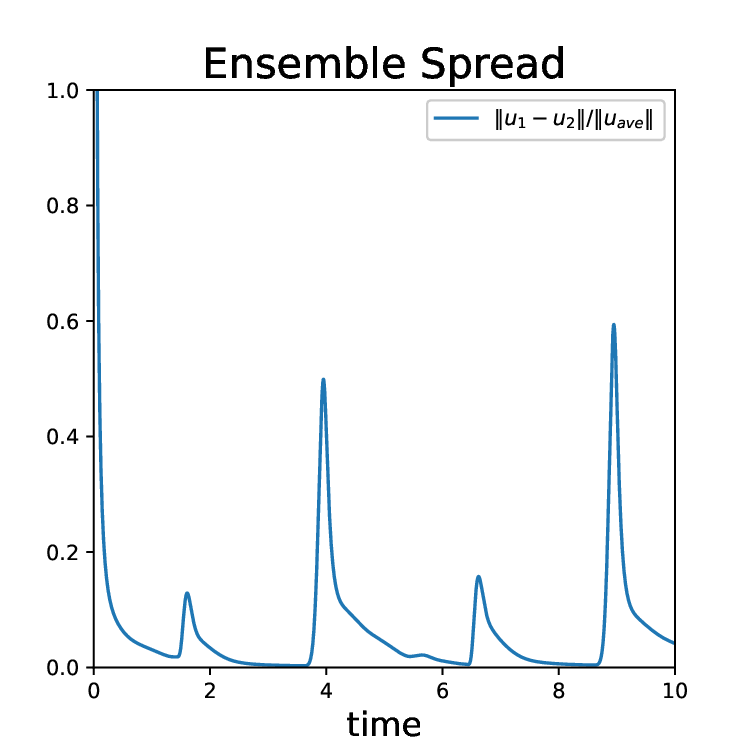}
   \caption{Ensemble spread.}
   \label{fig: victor-spread} 
\end{subfigure}
\hfill
\begin{subfigure}{0.45\linewidth}
   \centering
   \includegraphics[width=0.7\linewidth]{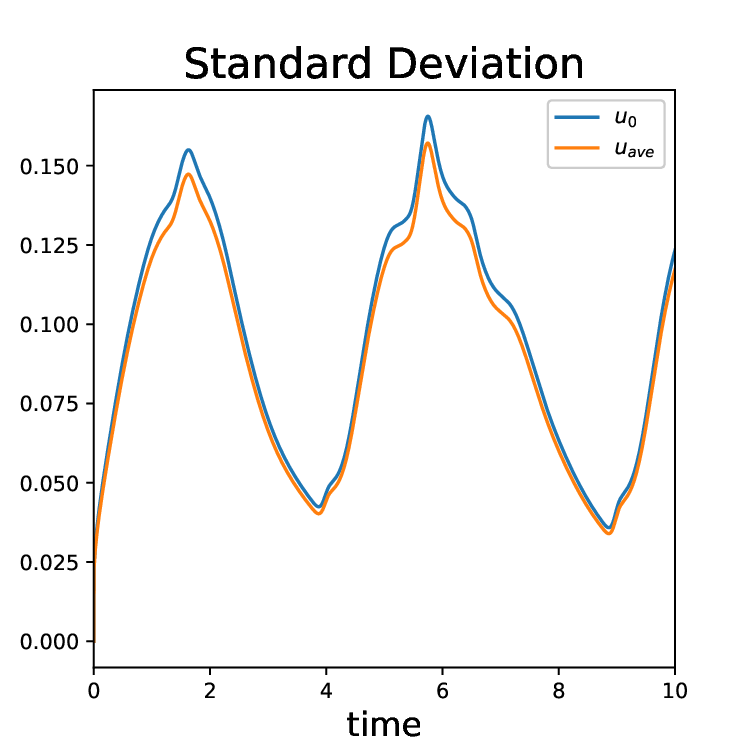}
   \caption{Standard deviation.}
   \label{fig: victor-std}
\end{subfigure}
\caption{}
\label{fig: spread-std-victor}
\end{figure}

Figure \ref{fig: victor-spread} shows that the ensemble spread changes periodically, with the peak of the spread approximately at $0.6$. We calculate the standard deviations considering \( u_0 \) as the mean and the ensemble mean \( u_{ave} \). Figure \ref{fig: victor-std} shows that the standard deviations for $u_0$ and $u_{ave}$ are similar. It indicates that the velocity is not chaotic.

\begin{figure}
\centering
   \begin{subfigure}{0.45\linewidth}
   \centering
\includegraphics[width=0.7\linewidth]{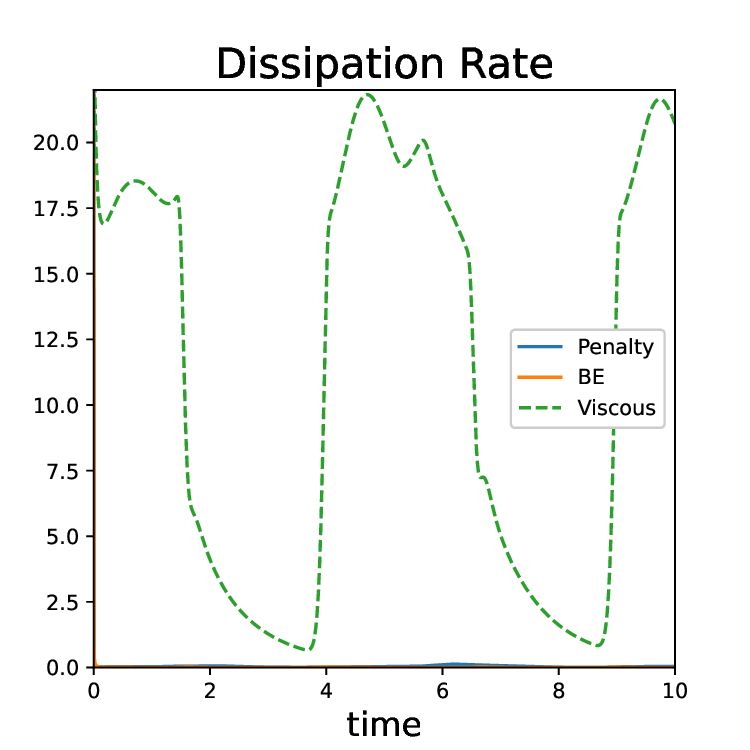}
   \caption{Dissipation rates.}
   \label{fig: aref-dissipation} 
\end{subfigure}
\hfill
\begin{subfigure}{0.45\linewidth}
   \centering
   \includegraphics[width=0.7\linewidth]{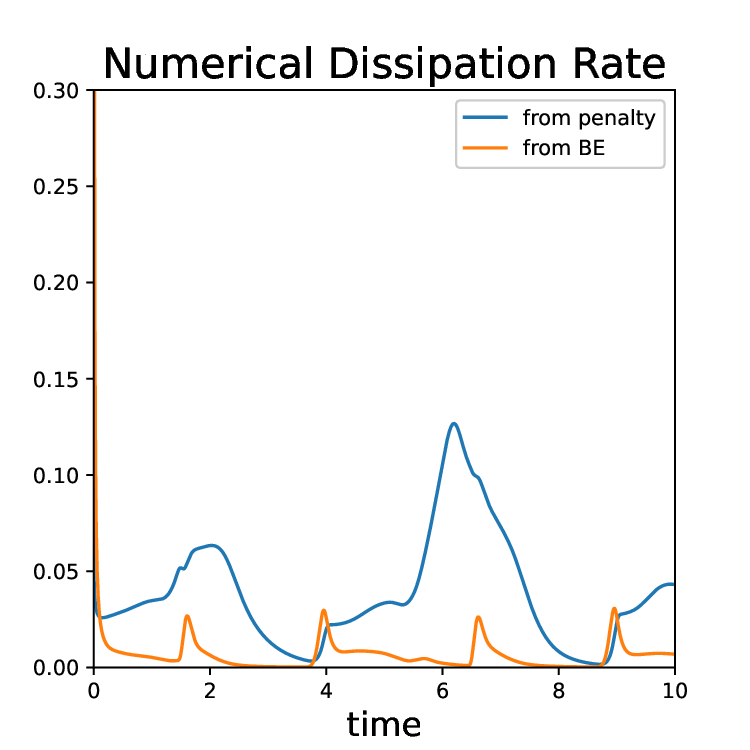}
   \caption{Numerical dissipation rates.}
   \label{fig: numerical_dissip}
\end{subfigure}
\caption{}
\end{figure}

In Figure \ref{fig: aref-dissipation}, we plot the numerical dissipation rates caused by penalizing the incompressibility condition and the BE time discretization. We compare them with the viscous dissipation rate. The numerical dissipation rate is much smaller than the viscous dissipation rate. In Figure 
\ref{fig: numerical_dissip},  the numerical dissipation rates have similar magnitudes and vary over time.

\begin{figure}
\centering
  \begin{subfigure}{0.45\linewidth}
   \centering
\includegraphics[width=0.7\linewidth]{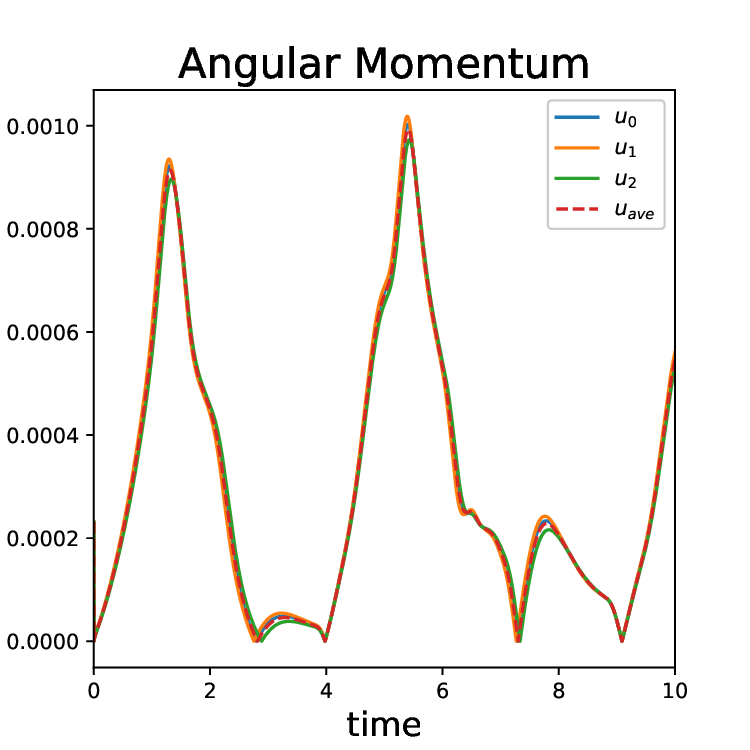}
   \label{fig:aref-angular-momentum}
   \caption{Angular momentum.}
   \end{subfigure}
\hfill
  \begin{subfigure}{0.45\linewidth}
   \centering
   \includegraphics[width=0.7\linewidth]{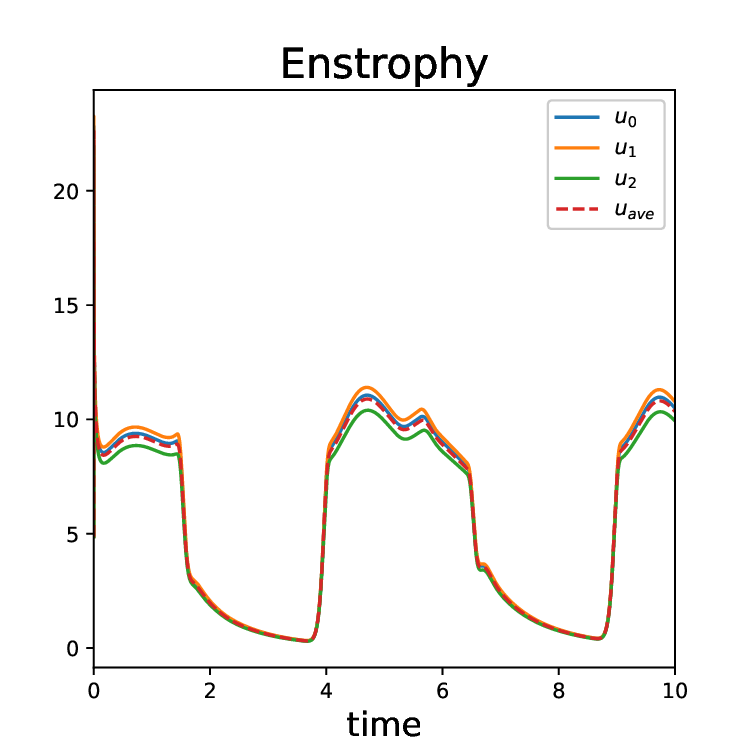}
   \label{fig:aref-enstrophy}
   \caption{Enstrophy.}
   \end{subfigure}
\vskip\baselineskip
\begin{subfigure}{0.45\linewidth}
   \centering   \includegraphics[width=0.7\linewidth]{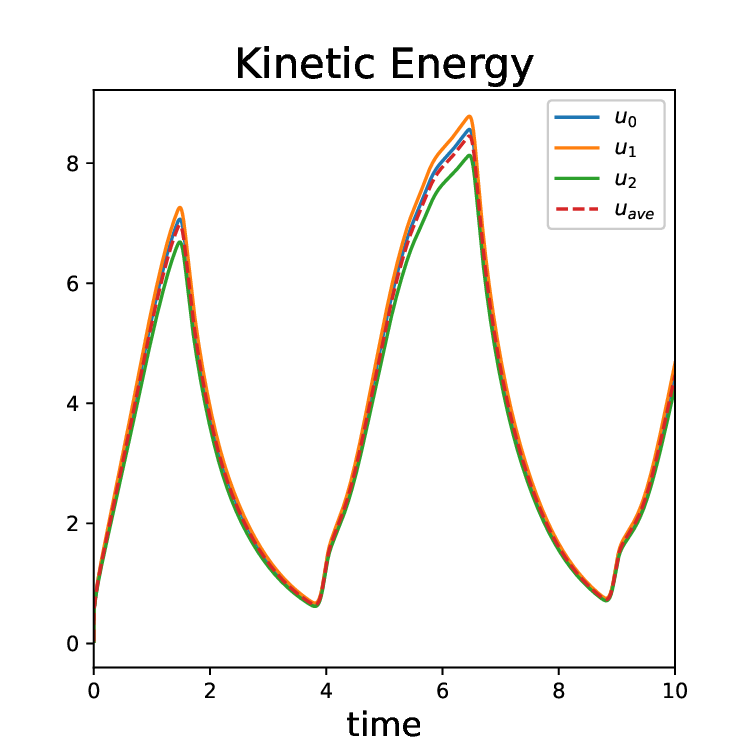}
   \label{fig:aref-kinentic-engery} 
   \caption{Kinetic energy.}
   \end{subfigure}
\hfill
  \begin{subfigure}{0.45\linewidth}
   \centering
\includegraphics[width=0.7\linewidth]{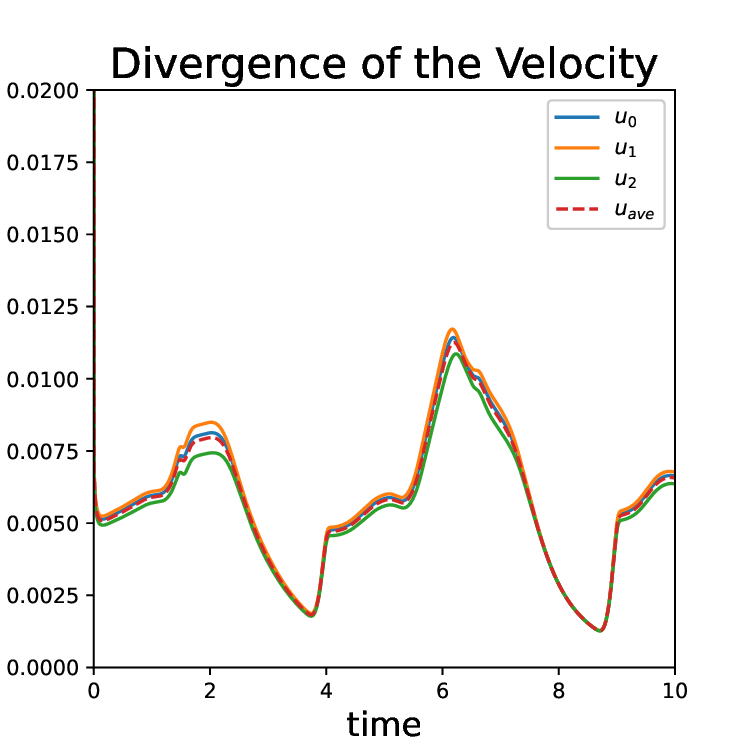}
   \label{fig:aref-divergence}
   \caption{Divergence of the velocity.}
   \end{subfigure}
\caption{Flow statistics for $u_{0}$, $u_{1}$, $u_{2}$ and $u_0$.}
\label{fig: victor-flow-statistics}
\end{figure}

We observe changes in kinetic energy, velocity divergence, angular momentum, and enstrophy as we activate and deactivate the left and right circles over time. The flow statistics of \( u_0 \), \( u_1 \), \( u_2 \), and \( u_{ave} \) are closely aligned in Figure  \ref{fig: victor-flow-statistics} and indicate the velocity field is not chaotic, where the trajectories of fluid particles exhibit chaotic behavior.
%\FloatBarrier

\subsection{Flow past a cylinder with the Coriolis force for large ensemble sizes}
\label{subsec: test-3}

Everything on Earth is rotating even without our noticing. The rotation changes the airflow and affects the climate, as discussed in Lee, Ryi, and Lim \cite{lee2017solutions}. The NSE with the Coriolis force is defined as follows:
\begin{equation*}
\frac{\partial u}{\partial t} + u \cdot \nabla u -\nu \Delta u +\nabla p +\omega Q u = f,
\end{equation*}
where $Q$ is a skew--symmetric matrix with a matrix norm equal to one, and $\omega$ is the Coriolis coefficient.

We extend the penalty--based ensemble method to the NSE with the Coriolis force. We evaluate this method using the benchmark 2D test flow past a cylinder, as described in \cite{schafer1996benchmark}. The inlet flow velocity is 
\begin{equation*}
u(x,y,t) = \left(\frac{6 y(0.41-y)}{0.41^2}, 0\right)^\top.
\end{equation*}

\begin{figure}
\centering
   \begin{subfigure}{0.45\linewidth}
   \centering
\includegraphics[width=0.7\linewidth]{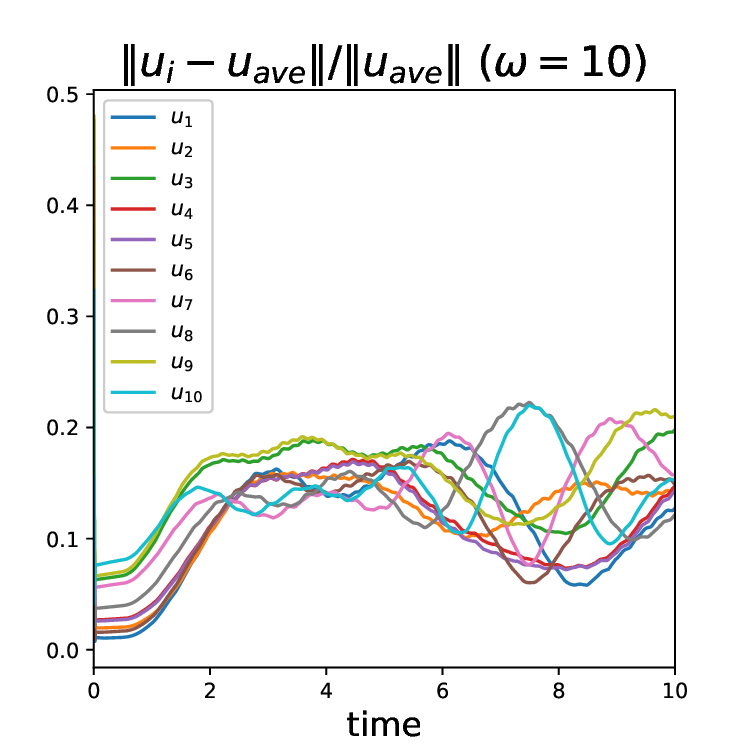}
   \caption{Spaghetti plot of the relative error.}
   \label{fig: spread} 
\end{subfigure}
\hfill
\begin{subfigure}{0.45\linewidth}
   \centering
   \includegraphics[width=0.7\linewidth]{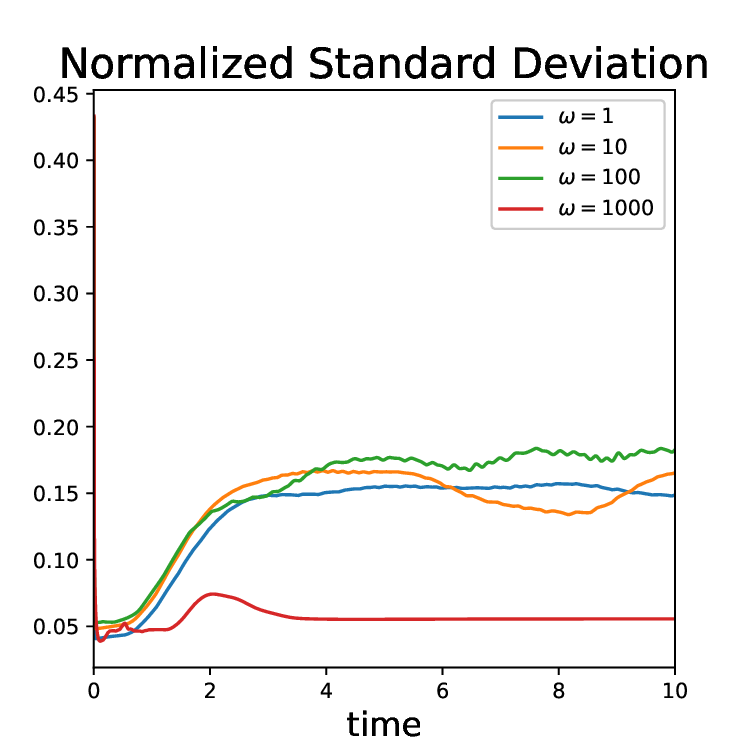}
   \caption{Normalized standard deviation.}
   \label{fig: std}
\end{subfigure}
\caption{The normalized standard deviation of the ensembles for different Coriolis coefficients.}
\label{fig: std-flow-past}
\end{figure}

\begin{figure}
\centering
   \begin{subfigure}{0.45\linewidth}
   \centering
   \includegraphics[width=0.7\linewidth]{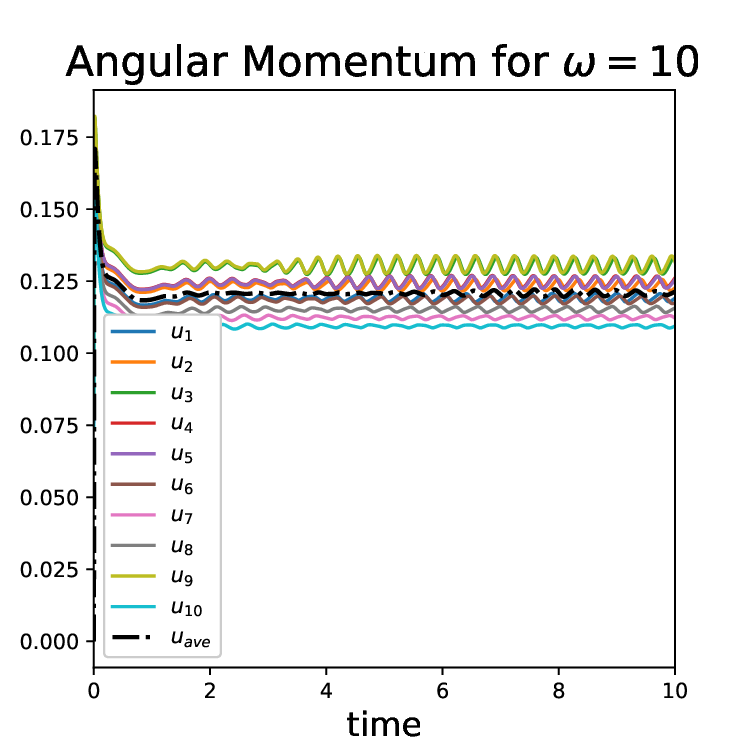}
   \caption{Angular momentum.}
   \label{fig:angular-momentum}
   \end{subfigure}
\hfill
   \begin{subfigure}{0.45\linewidth}
   \centering
   \includegraphics[width=0.7\linewidth]{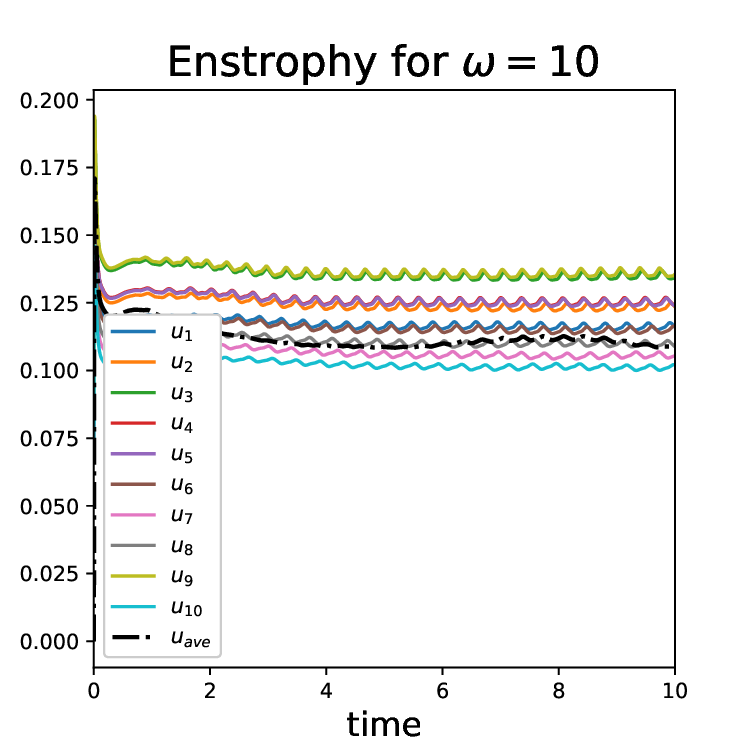}
   \caption{Enstrophy.}
   \label{fig:enstrophy}
   \end{subfigure}
\vskip\baselineskip
   \begin{subfigure}{0.45\linewidth}
   \centering
   \includegraphics[width=0.7\linewidth]{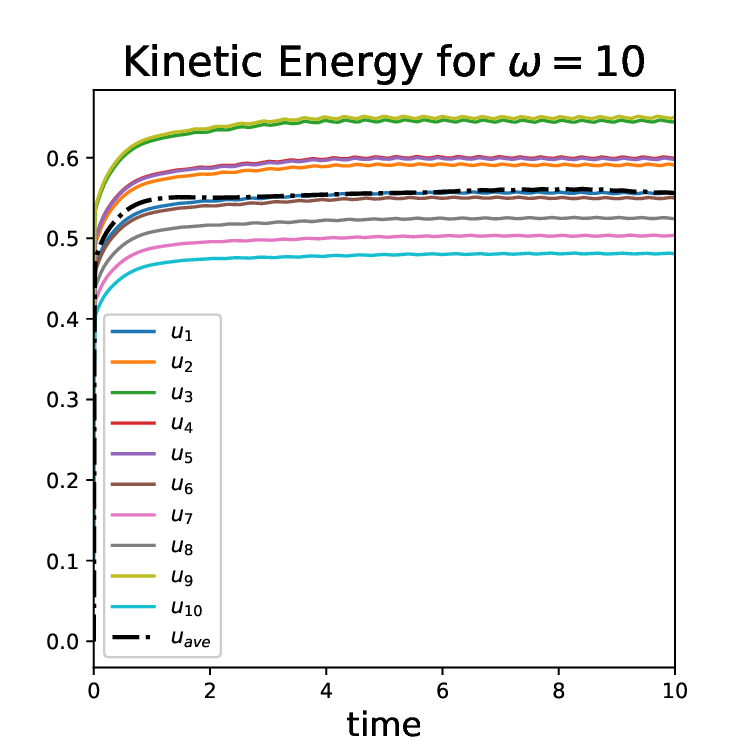}
   \caption{Kinetic energy.}
   \label{fig:energy}
   \end{subfigure}
\hfill
   \begin{subfigure}{0.45\linewidth}
   \centering
   \includegraphics[width=0.7\linewidth]{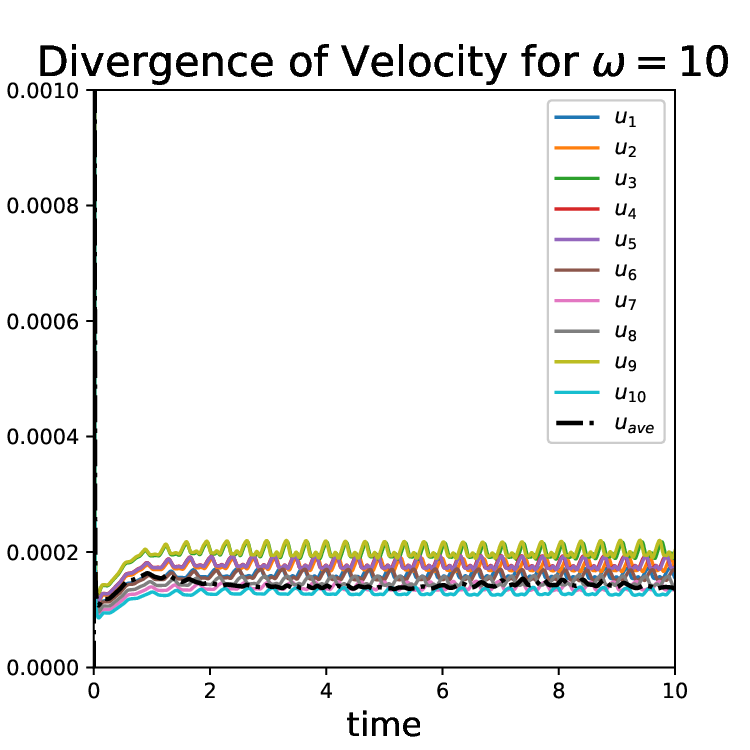}
   \caption{Divergence of the velocity.}
   \label{fig:div}
   \end{subfigure}
\caption{Flow statistics for all ensemble members and the mean flow at $\omega=10$.}
\label{fig: flow-past-all-members}
\end{figure}

We applied no--slip boundary conditions at the walls and on the obstacle. We generated second--order quadrilateral elements. We chose $J=10$, $T=10$, $\Delta t = 0.002$, $\nu = 0.001$, and $\epsilon = \Delta t$. The flow was at rest at $t=0$. We perturbed the inlet flow velocity for ensemble members as follows:
\begin{equation*}
u_j(x,y,t) = (1+\sigma_j \sin(2\pi y)) u, \text{ where } \ j=1, \ldots, 10.
\end{equation*}
$\sigma_j$ was randomly sampled from $-0.1$ to $0.1$. We first set $\omega =10$. Figure \ref{fig: spread} shows the spaghetti plot of the relative error of each single ensemble member to the mean flow. The normalized standard deviation for $\omega =10$ is around $0.15$ after $t=2$, as shown in Figure \ref{fig: std}. We calculated the angular momentum, enstrophy, kinetic energy, and velocity divergence for all ensemble members and the mean flow, as shown in Figure \ref{fig: flow-past-all-members}.

\begin{figure}
\centering
   \begin{subfigure}{0.45\linewidth}
   \centering
   \includegraphics[width=0.7\linewidth]{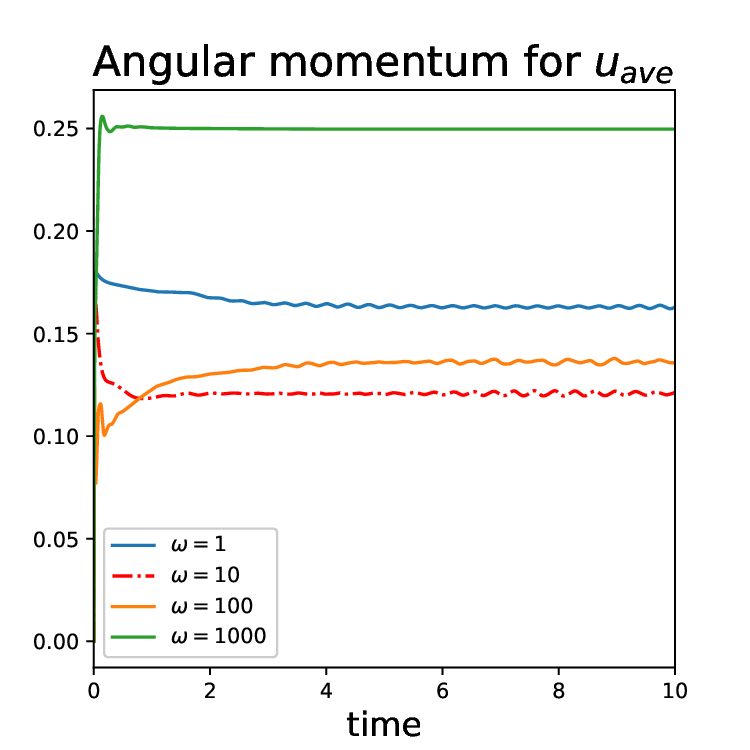}
   \caption{Angular momentum.}
   \label{fig:all-angular}
   \end{subfigure}
\hfill
   \begin{subfigure}{0.45\linewidth}
   \centering
   \includegraphics[width=0.7\linewidth]{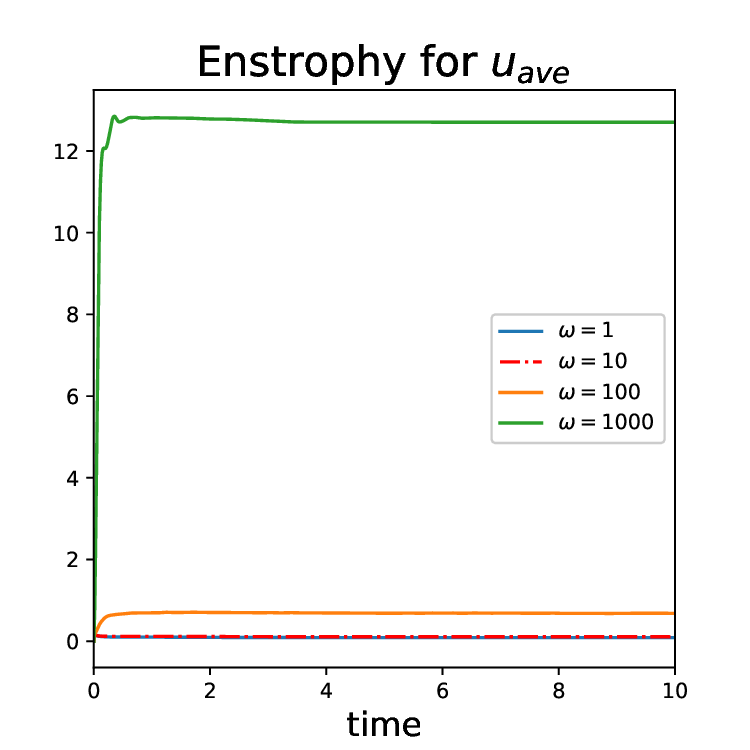}
   \caption{Enstrophy.}
   \label{fig:all-enstrophy}
   \end{subfigure}
\vskip\baselineskip
   \begin{subfigure}{0.45\linewidth}
   \centering
   \includegraphics[width=0.7\linewidth]{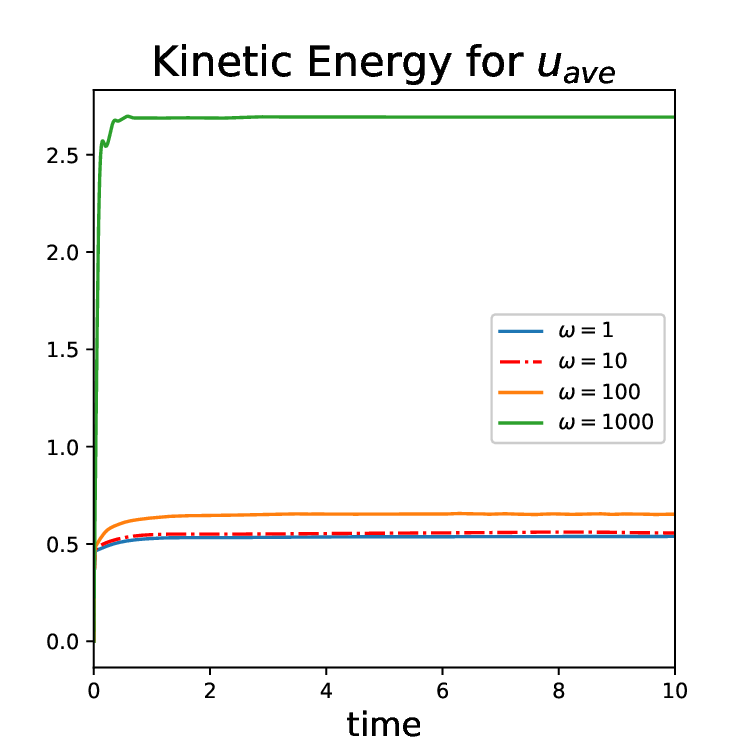}
   \caption{Energy.}
   \label{fig:all-energy}
   \end{subfigure}
    \hfill
   \begin{subfigure}{0.45\linewidth}
   \centering
   \includegraphics[width=0.7\linewidth]{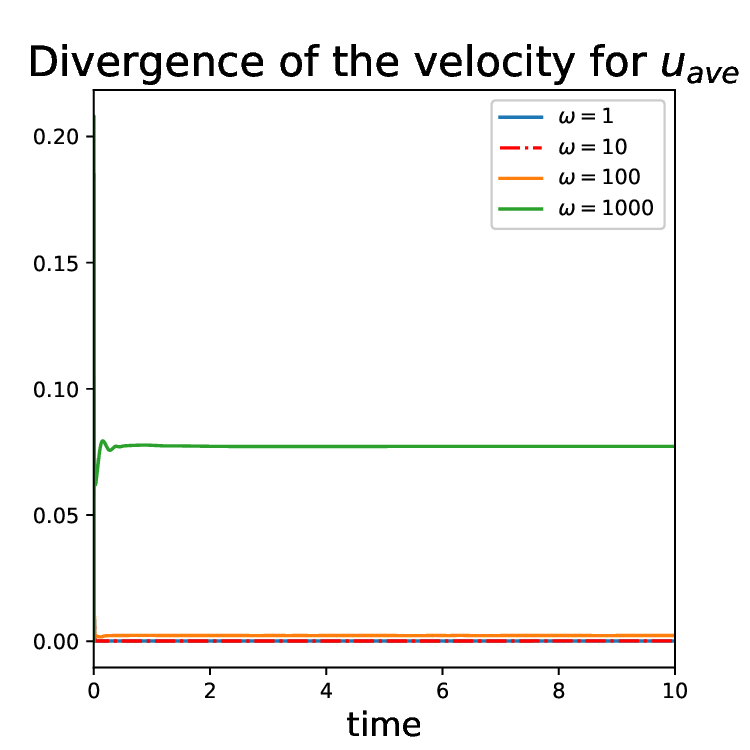}
   \caption{Divergence of the velocity.}
   \label{fig:all-div}
   \end{subfigure}
\caption{Flow statistics for the ensemble mean with different Coriolis coefficients.}
\label{fig: statistic-flow-past-cylinder}
\end{figure}

We set the Coriolis coefficient $\omega = 1, 10, 100,$ and $1000$ to study the effect of the Coriolis force. We calculate the normalized standard deviation for different values of the Coriolis coefficient, as shown in Figure \ref{fig: std}. For smaller $\omega$ values ($\omega = 1, 10,$ and $100$), the standard deviations are similar, around $0.15$. When increasing $\omega$ to $1000$, the rotational force becomes significant, resulting in a much smaller standard deviation. This indicates that the flow behaves like rigid body rotation. Additionally, we observe much larger magnitudes of angular momentum, enstrophy, kinetic energy, and divergence of the velocity for the ensemble mean when $\omega = 1000$, as shown in Figure \ref{fig: statistic-flow-past-cylinder}.

\section{Conclusions and prospects} 
Due to the chaotic nature of turbulent flows, numerical models have a finite predictability horizon. This predictability relies on the accuracy of the initial conditions. Small imperfections in the initial conditions can lead to losing predictive skill. While ensemble methods effectively address this issue, they can be computationally costly. To reduce the computational cost of ensembles while preserving accuracy. This method uses a shared coefficient matrix for all ensemble members. And it relaxes the incompressibility condition, uncoupling the flow velocity and pressure, thereby reducing model complexity, and allowing for a larger ensemble size.

We presented the stability and error estimates of the penalty--based ensemble method in equation (\ref{weak_form_method}). We extend the method to the NSE with random body forces and initial conditions with Monte Carlo sampling in Section \ref{sec: monte-carlo}. In Section \ref{subsec:2}, we verified the convergence rates with numerical experiments. In addition, we conducted a numerical experiment on chaotic advection, where the trajectories of the flow particles are chaotic, in Section \ref{sec: victor-aref}. Furthermore, we performed a benchmark test for flow past a cylinder with the Coriolis force using large ensemble sizes in Section \ref{subsec: test-3}.

Open problems include extending the penalty--based ensemble method to turbulence models with a higher Reynolds number \cite{fang20231, han2024numerical} and adapting penalty parameters for penalty--based ensemble methods.

\section*{Acknowledgments}
I thank my advisor, Professor William Layton, for his guidance and support. We thank Victor DeCaria for a helpful discussion of the test in Section \ref{sec: victor-aref}. The NSF partly supported this research of the author under grants DMS 2110379 and 2410893.

\bibliography{mybib}

\end{sloppypar}

\end{document}